\documentclass[ejs,noshowframe,preprint]{imsart}

\RequirePackage{amsthm,amsmath,amsfonts,amssymb}
\RequirePackage[numbers]{natbib}
\RequirePackage[colorlinks,citecolor=blue,urlcolor=blue]{hyperref}
\RequirePackage{graphicx}

\usepackage[utf8]{inputenc}
\usepackage{amsmath}
\usepackage{amsfonts}
\usepackage{amssymb}
\usepackage{amsthm}
\usepackage{graphicx}
\usepackage{float}
\usepackage{hyperref}
\usepackage[english]{babel}
\usepackage{enumitem}

\startlocaldefs

\renewcommand{\P}{\mathbb{P}}
\newcommand{\E}{\mathbb{E}}
\newcommand{\R}{\mathbb{R}}
\newcommand{\N}{\mathbb{N}}
\newcommand{\M}{ \mathbb{M}}
\newcommand{\m}{ \mathbf{m}}
\newcommand{\X}{\mathbf{X}}
\newcommand{\Z}{\mathbf{Z}}
\newcommand{\V}{\mathbf{V}}
\newcommand{\vv}{\mathbf{v}}
\newcommand{\U}{\mathbf{U}}
\newcommand{\x}{\mathbf{x}}
\newcommand{\z}{\mathbf{z}}
\newcommand{\y}{\mathbf{y}}
\newcommand{\Y}{ \mathbf{Y}}
\newcommand{\1}{\mathbf{1}}
\newcommand{\ck}{\check}

\renewcommand{\L}{\mathcal{L}} 

\newcommand{\Xx}{\underline{\X}}
\newcommand{\xx}{\underline{\x}}
\newcounter{exno}
\newcommand{\ex}[2]{\refstepcounter{exno}\label{#1} \noindent  \textbf{Example \ref{#1} (#2).\ }}
\newcommand{\rex}[1]{ \noindent  \textbf{Example \ref{#1} (continued).\ }}

\newtheorem{theorem}{Theorem}
\newtheorem{prop}{Proposition}
\newtheorem{lem}{Lemma}

\newtheorem{cor}{Corollary}
\newtheorem{rem}{Remark}

\endlocaldefs

\begin{document}

\begin{frontmatter}
\title{Parametric estimation and LAN property of the birth-death-move process with mutations}
\runtitle{Parametric estimation of the BDMM process}

\begin{aug}
\author[A]{\fnms{Lisa}~\snm{Balsollier}\ead[label=e1]{lisa.balsollier@cnrs.fr}}
\and
\author[A,B]{\fnms{Fr\'ed\'eric}~\snm{Lavancier}\ead[label=e2]{frederic.lavancier@ensai.fr}
}

\address[A]{Nantes Universit\'e, UMR 6629 CNRS\printead[presep={,\ }]{e1}}

\address[B]{Univ Rennes, Ensai, CNRS, CREST\printead[presep={,\ }]{e2}}
\runauthor{L. Balsollier and F. Lavancier}
\end{aug}

\begin{abstract}
A birth-death-move process with mutations is a Markov model for a system of marked particles in interaction, that move over time, with births and deaths. In addition the mark of each particle may also change, which constitutes a mutation. Assuming a parametric form of the model, we derive its likelihood expression and establish its local asymptotic normality. The efficiency and asymptotic distribution of the maximum likelihood estimator, with an explicit expression of its covariance matrix, are deduced. 
The underlying technical assumptions are shown to be satisfied by several natural parametric specifications. As an application, we leverage this model to analyse the joint dynamics of two types of proteins in a living cell, that are involved in the exocytosis process. Our approach enables to quantify the so-called colocalization phenomenon, answering an important question in cellular biology. 

\end{abstract}



\end{frontmatter}

\section{Introduction}

We address the parametric inference of birth-death-move processes with mutations, in short BDMM. These processes model the dynamics of a system of particles in $\R^d$ that move over time, while some new particles may appear and existing ones may disappear. In addition, each particle may be marked by a  label, that can change over time, which we call a mutation. This type of dynamics is observed in 
epidemiology \cite{masuda2017},  in ecology \cite{renshaw2001, pommerening2019} and in  cellular biology \cite{ronanfrederic, balsollier_2023}.  For example, the application we will consider later concerns the dynamics of proteins in a living cell, involved in the exocytosis process and observed near the plasma membrane of the cell: these proteins move within the cell, and for biological and photochemical reasons, some of them disappear while others appear. Moreover the motion of each protein can be of three types, the regime of which constitutes its label, and this regime may change over time, in line with what we call a mutation. 

Formally, a BDMM process is characterised by several components. First, the intensity of births, of deaths, and of  mutations, that are functions of the particle configurations, determine the waiting time until the next birth (or death, or mutation, respectively). Second, transition kernels specify how a birth (or a death or a mutation) occurs when such an event happens. Third, a continuous Markov process governs the motion of the system of particles between successive jumps, whether these jumps correspond to births, deaths, or  mutations. Assuming a parametric form for each of these characteristics, we are interested by their inference  given a single realisation of the process on a finite time interval. In this contribution, we address this question by maximum likelihood and we provide theoretical guarantees by proving the local asymptotic normality (LAN) property of the model, when the time interval increases and under some regularity conditions. This result yields the efficiency and asymptotic normality of the maximum likelihood estimator, with an explicit and estimable asymptotic covariance matrix.

Birth-death-move processes (without mutations) have been introduced in \cite{ronanfrederic} and further studied in \cite{emilien}, with a particular focus on their ergodic properties. These processes share close similarities with Markovian particle systems with killing and jumps as studied in \cite{eva2002, evalan},  branching processes \cite{athreya2012} and spatially structured population models  \cite{bansaye2015}. The introduction of possible mutations does not substantially change the  probabilistic properties of a birth-death-move process. In fact, a BDMM is a particular instance of the jump-move processes studied in \cite{emilien}. It follows that that it is Markovian whenever the intensity functions are bounded, and ergodic when the deaths compensate in some sense the births (precise conditions are stated in Appendix).

 From a statistical perspective,  no fully developed methodology currently exists for inference of birth-death-move processes (with or without mutations).  Nonparametric estimation of their intensity functions has been investigated in \cite{ronanfrederic}, but this contribution does not address the estimation of the transition kernels, nor of the inter-jump motions. On the other hand, parametric estimation has only been considered in the particular case when the process only involves spatial births and deaths, but no motions, and under some restrictive frameworks. Maximum likelihood inference in this setting and when $d=1$ has been considered in \cite{moller1994}, with an application to the analysis of dune displacements. The same  approach has been employed in \cite{Sadahiro19} for  a particular specification of the  process in $d=2$, aimed at  modeling openings and closures of shops in Tokyo. In these works, no theoretical study is provided. 
In \cite{evalan}, the LAN property of the closely related Markovian model of particle systems with killing and jumps is established, within a very general and abstract framework. In our  contribution, we derive the likelihood expression of the  general  birth-death-move process with mutations and we similarly show the LAN property under mild assumptions, which we show to be satisfied by standard parametric specifications of the model. 
As a result, we provide a comprehensive parametric approach to infer a general BDMM, with theoretical guarantees.  
 
The article is organised as follows. In Section~\ref{defBDMM}, the formal definition of the BDMM process is provided, together with specific examples for each of its characteristics.  In Section~\ref{sec:LAN}, we derive the likelihood expression of the model and we prove its LAN property. We in particular show that the associated technical hypotheses are satisfied for the natural examples presented in Section~\ref{defBDMM}, under mild assumptions. In Section~\ref{sec:simus}, a simulation study is carried out, illustrating the performances of the maximum likelihood estimator of some parameters of the model, along with the estimation of the corresponding  confidence regions. The setting is close to that of the dataset analysed subsequently, providing reassurance regarding the reliability of the conclusions concerning it. This dataset is a video sequence showing Langerin proteins and Rab-11 proteins in a living cell, acquired by fluorescence microscopy \cite{Boulanger2014}. As detailed in Section~\ref{section:simus:appli}, their dynamics are consistent with the realisation of BDMM process. Leveraging this model and its estimation, we prove that Langerin proteins are colocalized with Rab-11 proteins, and we quantify this phenomenon, answering an important biological question. Finally, Section~\ref{discussion} discusses our contribution and several perspectives, while an appendix gathers the proofs of our theoretical statements and some ancillary technical results.   

The Python code for  implementing  the simulation study and processing the data is available in our online GitHub repository at \url{https://github.com/balsollier-lisa/Parametric-estimation-of-the-BDMM-process}. 

\section{Definition of the BDMM process and examples}\label{defBDMM}

As mentioned in the introduction, a BDMM process is a Markov process describing the dynamics of a system of particles evolving over time. These dynamics include births and deaths, which modify the system's cardinality, as well as possible mutations of the particles' labels (if there are any). We first introduce in Section~\ref{sec:state} the state space of the process and some related notation. We then  describe  in Section~\ref{dynamique} the algorithmic construction of the process, providing  clear insight into its dynamics. 
Section~\ref{model examples} subsequently details each characteristic of the process (the jump intensities, the transition kernels and  the inter-jump motions) and presents several natural examples. 

\subsection{State space}\label{sec:state}

We denote by $(\X_{t})_{t \geq 0}$ the BDMM process. 
At each time $t \geq 0$, $\X_t$ describes the configuration of a set of marked particles. 
Each particle $x\in \X_t$ reads $x=(z,m)$ where $z\in\Lambda\subset \R^{d}$ encodes the spatial location of the particle along with its possible continuous mark, and where $m\in\M$ is a discrete mark (or \emph{label}), $\M$ being the finite set of possible labels.  

At each time $t$, $\X_t$ represents the set of all alive particles, where the ordering does not matter and the cardinality may change over time. For this reason, for any $n \geq 1$, we introduce the natural projection 
\begin{equation}\label{def pin}
\Pi_n : (x_1,\dots,x_n) \in ( \Lambda \times \M)^n \mapsto \{ x_1,\dots,x_n \}
\end{equation}
that identifies two elements $(x_1,\dots,x_n)$ and $(y_1,\dots,y_n)$ of $( \Lambda \times \M)^n$ if there exists a permutation $\pi$ of $\{1,\dots,n\}$ such that
 $x_i=y_{\pi(i)}$ for any $1 \leq i \leq n$.  

We set for any $ n \geq 1$, $E_n = \Pi_n( ( \Lambda \times \M)^n)$. The state space of $\X_t$ for $t \geq 0$ is then $$E= \bigcup_{n \geq 0} E_n,$$
 where $E_0=\{\text{\O}\}$ consists of the empty configuration. 
 
 The cardinality of a configuration $\x\in E$ will be denoted by $n(\x)$, so that $\x\in E_{n(\x)}$, and we write $$\x=\{x_1,\dots,x_{n(\x)}\}
=\{(z_{1},m_{1}),\dots, (z_{n(\x)},m_{n(\x)})\},$$ where $z_i\in \Lambda$ and $m_i\in\M$. In the following we will sometimes assimilate  $( \Lambda \times \M)^n$ with $ \Lambda^n \times \M^n$ for convenience, so that we may also write $\x=\{(\z,\m)\}$ where $\z=(z_1,\dots,z_{n(\x)})$ and $\m=(m_1,\dots,m_{n(\x)})$. Similarly, we will write the configuration at time $t$ of the BDMM process  in either way 
\[
	\X_{t}=\{x_{1,t}, \dots, x_{n,t}\} = \{(z_{1,t},m_{1,t}) \dots, (z_{n,t},m_{n,t})\} \quad\text{or}\quad\X_t=\{( \Z_t, \m_t)\},
\]
where $n=n(\X_t)$. Moreover, for $\x\in E_n$ and $(z,m)\in\Lambda \times \M$, we  write $\x\cup (z,m)$ for $ \{x_{1},\dots,x_{n},(z,m)\}\in E_{n+1}$,  and 
  for $n\geq 1$ and $1\leq i\leq n$, we write $\x\backslash x_{i}$ for $\{x_{1},\dots,x_{i-1},x_{i+1},\dots,x_{n}\}\in E_{n-1}$.

 We equip $E$ with the Borel $\sigma$-algebra $\mathcal E$ and we endow $E$ with the distance $d_1$ defined for $\x$ and $\y$ in $E$ such that $n(\x) \leq n(\y)$ by
\begin{equation*} \label{defd1}
     d_1(\x,\y) = \frac{1}{n(\y)} \left ( \min_{\pi \in S_{n(\y)}} \sum_{i=1}^{n(\x)} (\|x_i-y_{\pi(i)}\| \wedge 1) + (n(\y)-n(\x)) \right ),
    \end{equation*}
with $d_1(\x,\textnormal{\O})=1$ and where $S_n$ denotes the set of permutations of $\{1,\dots,n\}$. This distance makes in particular the function $\x\mapsto n(\x)$ continuous on $E$, and further implies some nice topological properties for $(E,d_1)$, see  \cite{schuhmacher2008} and \cite{emilien} for details.

\subsection{Dynamics of a BDMM process}\label{dynamique}

The dynamics of a BDMM process alternates continuous motions and jumps. Inter-jumps motions concern the location and/or the continuous mark of all alive particles, that move in $\Lambda$. Jumps are of three types: either a birth (when a new particle appears), or a death (when an existing particle disappears), or a mutation (when an existing particle changes its label). 
Accordingly, the dynamics is based on the three following elements:

 \begin{enumerate} 
 \item A continuous homogeneous Markov process $(\Y_{t})_{t\geq0}$ on $E$ that drives the motion of all particles of $(\X_{t})_{t\geq0}$ between two jumps. This process will not affect the labels of the particles, neither their cardinality.  Given an initial condition $\Y_0=\x$, $(\Y_{t})_{t\geq0}$ will be defined as the solution of a stochastic differential equation, as described and exemplified in Section~\ref{exmove}. 
  \item Three continuous non-negative functions $ \beta$, $ \delta$ and $ \tau$ on $E$, that refer respectively to the birth, death and mutation intensities. These functions govern the waiting times between the jumps of $(\X_t)_{t \geq 0}$. Heuristically, the probability that a birth occurs in the interval $[t,t+dt]$ given that the particles are in the configuration $\X_t$ à time $t$ is $\beta(\X_t)dt$, and similarly for $\delta$ and $\tau$. We will denote  $\alpha=\beta+\delta+\tau$ the total jump intensity of the process, that  we assume to be bounded. 
 \item Three transition kernels $K_{ \beta}$, $K_{ \delta}$, and $K_{ \tau}$ from $E\times \mathcal E$ to $[0,1]$, that specify how a jump occurs when it happens. For example, given that a birth occurs in a configuration $\x$, $K_\beta(\x,A)$ is the probability that the new configuration $\x\cup (z,m)$ belongs to $A\in \mathcal E$, where $(z,m)$ denotes the new particle, and similarly for $K_\delta$ and $K_\tau$. 
  \end{enumerate}
Details and examples of these three basic characteristics of the process are provided in the next section. 

Let us specifically describe the algorithmic definition of the process. This iterative construction follows \cite{ronanfrederic} and  \cite{emilien}, and may serve as a simulation procedure, see \cite{ronanfrederic} for details.  Let $(\Y^{(i)}_{t} )_{t\geq0}, i \geq 0$, be a sequence of processes on $E$, identically distributed as $(\Y_{t})_{t\geq0}$, in the sense that whatever $\x\in E$,  $(\Y^{(i)}_{t})_{t\geq0}$ given that $\Y^{(i)}_0=\x$  has the same distribution as  $(\Y_{t})_{t\geq0}$ given that $\Y_0=\x$. 

Starting from the initial configuration $\X_{0}$ at time $T_0=0$, we iteratively build the process as follows.
\begin{enumerate}[label=\roman*)]
    \item Given $\X_0$, generate the $n(\X_0)$ continuous trajectories of $(\Y^{(0)}_t)_{t\geq 0}$.
      \item Given $\X_0$ and $(\Y_t^{(0)})_{t \geq 0}$,  generate the first inter-jump time $T_1-T_0$ according to the cumulative distribution function
    \begin{equation}\label{jump law}F_1(t)=1-\exp\left(-\displaystyle\int_0^t \alpha(\Y^0_u)d u\right),\quad t>0.\end{equation}
The process until time $T_1$ is given by the generated trajectories, i.e. 
     $$(\X_t)_{T_0\leq t< T_1}=(\Y^0_{t-T_0})_{T_0\leq t< T_1}.$$
     
   \item Given $T_1$ and $\X_{T_{1}^-}=\Y^0_{T_1-T_0}$,  generate the first jump:
   \begin{itemize}
   \item it is a birth with probability $\beta(\X_{T_{1}^-})/ \alpha(\X_{T_{1}^-})$, in which case $\X_{T_1}$ is generated according to $K_{ \beta}(\X_{T_1^-},.)$;
   \item it is a death with probability $\delta(\X_{T_{1}^-})/ \alpha(\X_{T_{1}^-})$, in which case $\X_{T_1}$  is generated according to $K_{ \delta}(\X_{T_1^-},.)$;
   \item it is a mutation with probability $\tau(\X_{T_{1}^-})/ \alpha(\X_{T_{1}^-})$, in which case $\X_{T_1}$  is generated according to $K_{ \tau}(\X_{T_1^-},.)$.
   
   \end{itemize}
      
   \item Return to step i) with $T_0\leftarrow T_1$ and $\X_0\leftarrow \X_{T_1}$ to generate the new trajectories $(\Y^{(1)}_t)_{t\geq 0}$ starting from $\X_{T_1}$, the next jumping time $T_2$, and so on.
\end{enumerate}

The sequence of jumping times of the BDMM process $(\X_t)_{t \geq 0}$ is $(T_n)_{n \geq 1}$. 
For $t>0$, we denote by  $N_t$ the number of jumps on $[0,t]$, i.e. \[ N_t=\sum_{i\geq 1} \1_{ T_i\leq t}. \]
Similarly, we denote by $N_t^\beta$, $N_t^\delta$ and $N_t^\tau$ the number of births, deaths and mutations, respectively, before $t$. By our assumption that the total jump intensity $\alpha$ is bounded, we have that $N_t<\infty$ for any $t\geq 0$, i.e. there is no explosion of the process. 

The BDMM process, as defined by the above iterative construction, is a particular case of a jump-move process, as studied in \cite{emilien}. We deduce that it is a homogeneous Markov process with respect to its natural filtration $(\mathcal F_t)_{t\geq 0}$. We refer to the latter reference for the expression of the infinitesimal generator of the process and further probabilistic  properties.  
In the following, we will denote by $\P_\x$ and $\E_\x$ all probabilities and expectations given that $\X_0=\x$.

\subsection{Elements of the model and examples}\label{model examples}

In this section, we describe more precisely each characteristic of the process, namely the intensity functions, the transition kernels, and the motion process $(\Y_t)_{t \geq 0}$. 
We also present several illustrative examples, some of which were used in \cite{balsollier_2023} to successfully model  the dynamics of proteins inside a living cell (see also the application in Section~\ref{sec:simus}).

\subsubsection{Intensity functions}\label{ex intensity}

The birth, death and mutation intensity functions are such that, given $\X_{t}=\x$,  a birth (resp. a death or a mutation) occurs in $(t, t+h]$ with probability  $\beta(\x)h+o(h)$ (resp. $\delta(\x)h+o(h)$ or  $\tau(\x)h+o(h)$) as  $h\to 0$.
These intensities can also be seen in the following way: $ \beta(\X_{t^-})$ (resp. $ \delta(\X_{t^-})$ or $\tau(\X_{t^-})$) is the intensity of the associated counting process $N^ \beta_t$ (resp. $N^ \delta_t$ or $N^ \tau_t$). These interpretations are consequences of the specific form \eqref{jump law} of the waiting time before the next jump, see \cite{ronanfrederic} for details. 

\medskip

To set some examples, let us denote by $\gamma$ any of $\beta$, $\delta$ or $\tau$. The most simple situation is when the intensity function is a constant rate, that is $\gamma(\x)=\gamma$ for any $\x\in E$, where $\gamma>0$.  Then there is in average $\gamma\Delta$ new events that occur in any time interval of size $\Delta$, whatever the configuration of $\X_t$ is. Another typical setting is when the intensity is proportional to the cardinality, that is $\gamma(\x)=\gamma n(\x)$ for $\x\in E$ and $\gamma>0$. This reflects the situation where each particle has its own rate $\gamma>0$ and the particles do not interact, so that the total rate over all particles is $\gamma n(\x)$. These two examples are observed in the biological applications studied in \cite{ronanfrederic} and \cite{balsollier_2023}. More complicated examples of intensity functions are also considered in \cite{ronanfrederic}, where $\gamma(\x)$ depends on the underlying Voronoï tessellation induced by $\x$.

 \subsubsection{Transition kernels}\label{ex kernels}
 
Denoting $\gamma$ for any of $\beta$, $\delta$ or $\tau$, we assume that for any $\x$, the transition kernel $K_\gamma(\x,.)$ admits a density $k_\gamma(\x,\y)$, $\y\in E$, with respect to a measure $\nu_\gamma(\x,.)$, that is for all $\x \in E$ and $\y\in E$, $$K_\gamma(\x,d\y)=k_\gamma(\x,\y)\nu_\gamma(\x,d\y).$$ 
We now specify $k_\gamma$ and $\nu_\gamma$ in each case, whether the transition concerns: (i) a birth, (ii) a mutation or (iii) a death. 

\medskip

(i) For a birth transition, we assume that for any $F\in\mathcal E$, $\nu_\beta(\x,F)=0$ except if there exist $A\subset \Lambda$ and $ \mathcal{I}\subset  \M$ such that $$F=\x\cup (A \times \mathcal{I}):=\left\{	\y\in E, \:\exists \:(z ,m )\in A\times \mathcal{I},\: \y=\x\cup (z,m)\right\}.$$ In this case we set $\nu_\beta(\x,\x\cup(A \times \mathcal{I}))=|A|\times |\mathcal I|$ to be the Lebesgue measure on $\Lambda\times \M$. This choice of $\nu_\beta$ ensures that the birth transition produces the addition of exactly one particle to the configuration $\x$, as expected for a birth.
For the density $k_\beta(\x,\y)$, we assume that
\begin{align}\label{birth kernel}
k_{ \beta}(\x, \y)= \left\{
    \begin{array}{ll}
 \displaystyle p_{m}k_{ \beta}^m(\x,z)  &\text{if }\exists (z,m)\in\Lambda\times \M,\   \y=\x\cup (z,m),\\
   0& \text{otherwise,} 
    \end{array}
   \right.
\end{align}
where $p_{m} \in [0,1]$ is the probability that the new particle has the label $m \in \M$, with $\sum_{m\in \M}p_{m}=1$,  and $k_{ \beta}^m(\x,.)$ is the density for the location of the new particle in $\Lambda$ given that its label is $m$. 

A simple example of birth transition is the uniform law for the label and the location, which corresponds to  $p_{m}=  1/ | \M|$ for all $m\in \M$ and  $k_{ \beta}^{ m}(\x,z)= 1/ | \Lambda|$, for all $z \in \Lambda$ and $m\in\M$. More sophisticated examples are provided below.

\medskip

\ex{mixture}{Birth kernel as a mixture of normal laws} 
In this example, a new particle is more likely to appear close to existing particles. The birth density is a mixture of isotropic normal distributions, centred at each existing particle, with deviation $\sigma>0$. Specifically,  for any configuration $\x=\{ (z_{1},m_{1}),\dots, (z_{n(\x)},m_{n(\x)})\}$ and any $m\in \M$,
\begin{align*}
k_{ \beta}^{ m}(\x,z)= \frac 1 {n(\x)} \sum_{i=1}^{n(\x)} \frac{ 1}{ o(z_{i},z)}\exp\left(-\frac{\|z- z_i\|^2}{2\sigma^2}\right) ,
\end{align*}
where $o(z_{i},z)=\int_{ \Lambda} \exp\left(-\frac{\|z- z_i\|^2}{2\sigma^2}\right)dz.$ We may easily  extend this example  by requiring that a particle with label  $m$ can only appear close to particles with labels $m'$. \\
 
 \ex{potentiel}{Birth kernel driven by a potential} In this example the birth kernel is given through a function $V:E \to \R$ by:
$$k_{ \beta}^{ m}(\x, z)= \frac{ 1}{ c^m(\x)}e^{-(V(\x\cup (z,m))-V(\x))}$$
where $c^m(\x)= \int_{ \Lambda}e^{-(V(\x\cup (z,m))-V(\x))}dz$. Given a configuration $\x$, a new particle is more likely to appear in the vicinity of points $z \in \Lambda$ that make $V(\x\cup (z,m))-V(\x)$ minimal. Following the statistical physics terminology, the function $V$ can be interpreted as a potential that the system tends to minimise  at each birth.  A typical instance, for $\x=\{ (z_{1},m_{1}),\dots, (z_{n(\x)},m_{n(\x)})\}$,  is $V(\x)=\sum_{i\neq j} \Phi_{m_i,m_j}(z_i-z_j)$ for some pair potential functions $\Phi_{m,m'}$, $m,m'\in\M$, see \cite{emilien} for some examples.

\bigskip

(ii) Concerning the mutation transition kernel, we set for any $\x=\{x_1,\dots,x_{n(\x)}\}$, where $x_i=(z_i,m_i)$, and any $F\in\mathcal  E$, 
$$\nu_\tau(\x,F)=\sum_{i=1}^{n(\x)}\sum_{m\in\M}  \1_{(\x \setminus (z_i,m_i))\cup (z_i,m)\in F},$$ 
so that a transition only concerns a change of label of an existing particle.  For the density $k_\tau$, we assume that 
\begin{align*}
k_{ \tau}(\x,\y)=\left\{
    \begin{array}{ll}
s(x_{i},\x) q_{m}(x_{i},\x)  &\text{if }\exists x_i\in\x, \exists m\in \M,  \  \y=(\x\backslash (z_{i},m_i))\cup (z_{i},m),\\
   0& \text{otherwise.} 
    \end{array}
   \right.
\end{align*}
In this expression, $s(x_{i},\x) \in [0,1]$ is the probability that the particle $x_{i}$ in $\x$ changes his label and $q_{m}(x_{i},\x)$ is the probability that, given that the particle $x_{i}\in\x$ mutates, its label changes from $m_{i}$ to $m$. We thus have $\sum_{i=1}^{n(\x)} s(x_i,\x)=1$ and $\sum_{m\in \M} q_{m}(x_{i},\x)=1$ with $q_{m_i}(x_{i},\x)=0$. 
An example is provided below.

\medskip

\ex{transitionmatrix}{Mutation kernel given by a transition matrix} A natural example consists in choosing the particle to be modified uniformly among all particles, and then changing its label according to a transition matrix with coefficients $(p_{m,m'})_{m,m'\in \M}$. This corresponds to the choices 
$$s(x_{i},\x)= \frac{ 1}{n(\x)} \quad \text{and }\quad q_{m}(x_{i},\x)= p_{m_{i},m},$$ 
where $p_{m,m}=0$ and for all $m\in \M$, $ \sum _{m' \in \M}p_{m,m'}=1$.

\bigskip

(iii) Finally, for the death transition kernel, we set for any $\x=\{x_1,\dots,x_{n(\x)}\}$ and any $F\in\mathcal E$, $\nu_{ \delta}(\x,F)=\sum_{i=1}^{n(\x)} \1_{\x\backslash x_{i}\in F}$, so that this transition only concerns the death of an existing particle. For the density, 
\begin{align*}
k_{ \delta}(\x,\y)= \left\{
    \begin{array}{ll}
\omega(x_{i},x)  &\text{if } \y=\x\backslash x_{i},\\
   0& \text{otherwise,} 
    \end{array}
   \right.
\end{align*}
where $\omega(x_{i},\x) \in [0,1]$ is the probability that the particle $x_{i}$ in $\x$ dies, with $\sum_{i=1}^{n(\x)} \omega(x_i,\x)=1$. A simple example is the uniform death where $\omega(x_i,\x)=1/n(\x)$ for any $i$.

\subsubsection{The inter-jump motion}\label{exmove}

Remember that between two jumps, the  BDMM process $(\X_t)_{t\geq0}$ has a constant cardinality $n=n(\X_t)$ and constant labels in $\M$. To specify the inter-jump motion of $(\X_t)_{t\geq0}$, it is then enough to define the dynamics of a process $(\Y ^{|n}_{t})_{t\geq0}$ on each subspace $E_{n}$ where 
$$\Y^{|n}_t=\{(z_{1,t},m_1),\cdots,(z_{n,t},m_n)\}.$$
Accordingly, the position and the continuous mark of each particle may move, that is $z_{i,t}\in\Lambda$ may move, but the label $m_i$ will remain constant. To define such dynamics, we first need to come back to a standard system of $n$ ordered (labelled)  particles in $(\Lambda\times\M)^n$, that is 
$$\tilde \Y^{|n}_t=((z_{1,t},m_1),\cdots,(z_{n,t},m_n)).$$
The full definition of the move process $(\Y_t)_{t\geq 0}$ on $E$ then follows the three steps:
\begin{enumerate}
\item Define the dynamics of $\tilde \Y^{|n}_t$ on $(\Lambda\times \M)^n$ thanks to a system of stochastic differential equations, like $(M^{|n})$ specified below, where the motion acts only in $\Lambda^n$.
\item Provided this system yields a solution whose distribution satisfies the permutation equivariance property (see below), deduce $\Y^{|n}_t = \Pi_n(\tilde\Y^{|n}_t)$ on $E_n$, where $\Pi_n$ is given by \eqref{def pin}.
\item Then define $(\Y_t)_{t\geq 0}$ on $E$ by $\Y_t=\sum_{n\geq 0} \Y^{|n}_t \1_{\{\Y_0\in E_n\}}$.
\end{enumerate}
This construction and the following facts are detailed in \cite{emilien}. The permutation equivariance property means that for any permutation $\pi$ of $\{1,\dots,n\}$, the law of
$\tilde\Y^{|n}_{t}=((z_{1,t},m_1),\cdots,(z_{n,t},m_n))$ given  $\tilde\Y^{|n}_{0}=((z_{\pi(1)},m_{\pi(1)}),\cdots,(z_{\pi(n)},m_{\pi(n)}))$ is the same as the distribution of 
$((z_{\pi(1),t},m_{\pi(1)}),\cdots,(z_{\pi(n),t},m_{\pi(n)}))$ given $\tilde\Y^{|n}_{0}=((z_{1},m_1),\cdots,(z_{n},m_n))$. This property says that if we rearrange the ordering of the coordinates of the initial state, the ordering of the solution is rearranged in the same way. Importantly, this means that we are free to choose any  initial ordering when defining $\tilde \Y^{|n}_0$ from $\Y^{|n}_0$, since it will lead to the same final point configuration $\Y^{|n}_{t}$, see (2.10) in \cite{emilien} for a mathematical argument. 
 Under this assumption, the Markov process $\Y$ is well-defined on $E$ and its transition kernel $Q_t^\Y$, defined for any bounded and measurable function $f$ on $E$ by $Q_t^\Y f(x)=\E[ f(\Y_t)| \Y_0=x]$, reads 
\begin{equation*}\label{decomp_QY} Q_t^\Y f(x) = \sum_{n \geq 0} Q_t^{\tilde \Y^{ |n}}(f \circ \Pi_n)((x_1,\dots, x_n))\1_{x \in E_n},\end{equation*}
where $Q_t^{\tilde \Y^{|n}}$ denotes  the  transition kernel of $\tilde \Y^{ |n}$ in  $(\Lambda\times \M)^n$. Moreover, $\Y$ is continuous in $(E,d_1)$ whenever  for any $n$, $\tilde \Y^{ |n}$ is continuous in $(\Lambda\times \M)^n$ for the usual Euclidean distance.

As a consequence of this construction, the definition of $(\Y_t)_{t\geq 0}$ boils down to the first step above, that is the specification of the system of SDE $ (M^{|n})$, provided the latter is permutation equivariant. We assume in this paper that $\tilde \Y^{|n}_t=(\Z_{t},\m)$, starting at $t=0$ from $(\z,\m)$ where $\z=(z_1,\dots,z_n)$, is the solution of 
\begin{align*}
    M^{|n}(\z,\m):\quad \left\{
    \begin{array}{ll}
    d z_{i,t}=b_{i,n}(\Z_{t}, \m)d t +\sigma_{i,n}( \Z_t, \m)d B_{i,t},& t\geq 0, \quad i=1,\dots,n,\\
    z_{i,0}=z_{i},\quad i=1,\dots,n.& 
    \end{array}
   \right.
\end{align*} 
Here the drift functions $b_{1,n},\dots,b_{n,n}$ take their values in $\R^{d}$, the diffusion coefficients $\sigma_{1,n},\dots,\sigma_{n,n}$ are $(d,d)$ invertible matrices  and $(B_{1,t})_{t\geq0},\dots,(B_{n,t})_{t\geq0}$ are independent standard Brownian motions on $\R^{d}$. We  assume that the functions $b_{i,n}$ and $\sigma_{i,n}$ are globally Lipschitzian, so that a strong solution exists. 
If $\Lambda\neq \R^{d}$, we further assume that edge conditions (reflective or periodic) are added to ensure that the solution stays in $(\Lambda\times \M)^n$, see for instance \cite{fattler2007}.

The permutation equivariance property is ensured if for any $i=1,\dots,n$ and any  $\pi\in S_n$, 
\begin{equation*}
\begin{cases}
b_{i,n}((z_{\pi(1)},m_{\pi(1)}),\cdots,(z_{\pi(n)},m_{\pi(n)}))=b_{\pi(i),n}((z_{1},m_1),\cdots,(z_{n},m_n)),\\
\sigma_{i,n}((z_{\pi(1)},m_{\pi(1)}),\cdots,(z_{\pi(n)},m_{\pi(n)}))=\sigma_{\pi(i),n}((z_{1},m_1),\cdots,(z_{n},m_n)).
\end{cases}
\end{equation*}
This is for instance the case if $b_{i,n}(\Z,\m)=b(z_i,m_i)$ and  $\sigma_{i,n}(\Z,\m)=\sigma(z_i,m_i)$ for some functions  $b$ and $\sigma$, a situation where there is no interaction between the particles, as considered in the microbiological applications of \cite{briane} and \cite{balsollier_2023}. An example that includes interactions through a potential function is provided in the following example.

\medskip

\ex{langevin}{Langevin diffusion}  In this example the motion of each particle depends on an interaction force with the other particles, driven by a pair potential function $\Phi_{m,m'}$, $m,m'\in\M$, as in Example~\ref{potentiel}.  Specifically, in this model $\sigma_i(\Z,\m)= \sigma_{m_{i}}$ and 
$$b_{i,n}(\Z,\m)=- \sum_{j\neq i}  \Phi_{m_i,m_j}(z_{i}-z_{j}).$$
For existence, each pair potential $\Phi_{m,m'}$ must be smooth enough. We refer to \cite{fattler2007} and \cite{emilien} for details and examples.

\bigskip

For later purposes, let us specify the likelihood of $\Y^{|n}_t$ constructed as above. Let us first rewrite $(M^{|n})$ so that it takes the form of a standard stochastic differential equation. Denote by $\bar\sigma_n$ the block diagonal matrix of size $(nd,nd)$ formed by the $\sigma_{i,n}$ matrices, and  by $\bar b_n$ the vector of size $nd$ formed by the concatenation of the $b_{i,n}$ vectors.
Then, by denoting $B_{t}$ the vector formed by the concatenation of the vectors $B_{i,t}$, we have
\begin{align}\label{EDSglob}
    M^{|n}(\z,\m):\quad \left\{
    \begin{array}{ll}
    d  \Z_{t}=\bar b_n(\Z_{t}, \m)d t +\bar \sigma_{n}( \Z_t, \m)d B_{t},& t\geq 0, \\
     \Z_{0}= \z. &
    \end{array}
   \right.
\end{align} 
It is well known, cf \cite{liptser_1,klebaner}, that the Radon-Nikodym density of the solution $(\Z_{t})_{t\geq 0}$ of $M^{|n}(\z,\m)$ with respect to the reference process $\U_t=\z + \int_0^t \bar \sigma_{n}( \U_s, \m)d B_{s}$, on the interval $[0,t]$ reads
\begin{multline}\label{Lmove}
L(\Z_{[0,t]}, \m)=\exp\bigg(\displaystyle\int_{0}^{t} \bar b_n(\Z_s,\m)^T a_{n}^{-1}( \Z_s, \m)d  \Z_s  \\
-\frac{1}{2}\int_{0}^{t}\bar b_n( \Z_{s}, \m)^T a_{n}^{-1}( \Z_{s}, \m) \bar b_n( \Z_s,\m)d s\bigg),
\end{multline} 
where $a_{n}(\z,\m)= \bar\sigma_{n}(\z,\m) \bar\sigma_{n}^T(\z,\m)$. Note that the reference process must depend on the diffusion coefficient $\sigma_{n}$, since in continuous time the laws of two processes with  different diffusion coefficients are mutually singular. From a statistical point of view, this is not a restriction, as $\sigma_n$ can be considered known under continuous time observations, see for instance  \cite[Remark 10.3 in Chapter 10.6]{klebaner}. From \eqref{Lmove}, we get the density of the process $\tilde \Y^{|n}_t=(\Z_{t},\m)$ on $[0,t]$ with respect to the reference process $(\U_t,\m)_{t\geq 0}$ that belongs to $(\Lambda\times \M)^n$. 
Finally, by our construction above, if  $\Y^{|n}_t=\Pi_n((\Z_t,\m))$, then the density of $\Y^{|n}_t$ with respect to $\Pi_n((\U_t,\m))$ on $[0,t]$ takes the same form and does not depend on the ordering chosen to define $\Z_t$, $\z$ and $\m$ from $\Y^{|n}_t$, thanks to the permutation equivariance property.

\section{Likelihood and LAN property}\label{sec:LAN}

\subsection{Likelihood of the BDMM process}

As described in the previous section, the BDMM process $(\X_t)_{t\geq 0}$ depends on the intensity functions $\beta$, $\delta$, $\tau$, assumed to be bounded on $E$, on transition kernel densities for the births, the deaths and the mutations, denoted by $k_\beta(\x,\y)$, $k_\delta(\x,\y)$  and $k_\tau(\x,\y)$, as detailed in Section~\ref{ex kernels}, and on a continuous Markov diffusion model on $E$ that drives the inter-jump motions, see Section~\ref{exmove}.

Remember that between two jumps $T_i$ and $T_{i+1}$, $\X_t=\Y^{(i)}_{t-T_i}$ where $\Y^{(i)}$ has the same law as $\Y$. Write $\Y^{(i)}_t=\{(\Z^{(i)}_t,\m^{(i)})\}$ where $\m^{(i)}$ does not depend on $t$, as required for the labels of the continuous Markov process that drives the inter-jump motions. As assumed in Section~\ref{exmove}, $\Z^{(i)}$ is the solution of the stochastic differential system $M^{|n}(\z^{(i)},\m^{(i)})$ given by \eqref{EDSglob}, where $(\z^{(i)},\m^{(i)})$ is the initial configuration, i.e. $\X_{T_i}=\{(\z^{(i)},\m^{(i)})\}$, and $n=n(\X_{T_i})$. The likelihood of the inter-jump motion between $T_i$ and $T_{i+1}$, given $T_i$,  $T_{i+1}$ and $\X_{T_i}$, is thus $L(\Z^{(i)}_{[0,T_{i+1}-T_{i}]}, \m^{(i)})$ given by \eqref{Lmove}. We will simply write the latter $L(\X_{[T_i,T_{i+1}[})$  in the following, since $\X_t=\{(\Z^{(i)}_{t-T_i},\m^{(i)})\}$ for all $t\in[T_i,T_{i+1}[$.

Assume that the BDMM is parametrized by $\vartheta\in\Theta\subset \R^\ell$, for a given $\ell>0$. This means that all features of the BDMM may possibly depend on $\vartheta$, to the exception of the diffusion coefficient $\bar\sigma_n$ in \eqref{EDSglob} and \eqref{Lmove}. As stressed in the previous section, this is indeed a common fact that for continuous time observations, the diffusion coefficient can be considered  known,  so that it makes no sense to conduct inference on it in this setting (the story is of course different for discrete time observations). 
So we assume that the intensities, the kernel densities and the drift coefficients may depend on $\vartheta$.  
To emphasize this dependence, we introduce the argument $\vartheta$ into the notation of each of them.
 For example we will write $\beta(\vartheta,\x)$ instead of $\beta(\x)$, $k_\beta(\vartheta,\x,\y)$ instead of $k_\beta(\x,\y)$, $\bar b_n( \vartheta,\x)$ instead of  $\bar b_n(\x)$, $L(\vartheta,\X_{[T_i,T_{i+1}[})$ instead of $L(\X_{[T_i,T_{i+1}[})$, and so on.

%

Let us denote the set of birth times up to time $t$  by $$\mathcal{B}_{t}=\left\{ i \in \N, \: T_{i}\leq t \text{ and }T_i\text{ is a birth time}\right\},$$ and similarly $\mathcal{D}_{t}$ and $\mathcal{T}_{t}$  for the set of death times and mutation times up to time $t$. 

The following theorem provides the likelihood of the BDMM process. When there is no move and no mutation, we recover the formula given in  \cite{moller1994} for their applications to the dynamics of dunes in dimension $d=1$. It also has a consistent expression with the likelihood of a system of particles with killing and jumps, as derived in \cite{eva2002}. The justifications, detailed in Section~\ref{appendix:likelihood}, are however different and simpler, as they only rely on successive backward applications of  conditional expectations.  
\begin{theorem}\label{densite}
    Let $(\X_t)_{t\geq0}$ be a BDMM parameterized by $\vartheta \in\Theta$ and observed in continuous time on $[0,t]$, for $t>0$.  Then, the  likelihood $ \mathcal{L}_{t}(\vartheta)$ is expressed as follows:
\begin{align*}
   \mathcal{L}_{t}(\vartheta)&=\mathcal{L}^{\mathcal{B}}_{t}(\vartheta)\mathcal{L}^{\mathcal{D}}_{t}(\vartheta)\mathcal{L}^{\mathcal{T}}_{t}(\vartheta) \mathcal{L}_{t} ^{\text{move}}( \vartheta),
       \end{align*}
with
        \begin{align*}
      &\mathcal{L}^{\mathcal{B}}_{t}(\vartheta)=   \exp\left(-\displaystyle\int_0^{t}\beta(\vartheta,\X_{s})d s\right)\prod_{i\in \mathcal{B}_{t}}k_{ \beta}(\vartheta,\X_{T_{i^-}},\X_{T_{i}})\beta(\vartheta,\X_{T_{i^-}}),\\
       &\mathcal{L}^{\mathcal{D}}_{t}(\vartheta)=   \exp\left(-\displaystyle\int_0^{t}\delta(\vartheta,\X_{s})d s\right)\prod_{i\in \mathcal{D}_{t}}k_{ \delta}(\vartheta,\X_{T_{i^-}},\X_{T_{i}})\delta(\vartheta,\X_{T_{i^-}}),\\
       &\mathcal{L}^{\mathcal{T}}_{t}(\vartheta)=   \exp\left(-\displaystyle\int_0^{t}\tau(\vartheta,\X_{s})d s\right)\prod_{i\in \mathcal{T}_{t}}k_{ \tau}(\vartheta,\X_{T_{i^-}},\X_{T_{i}})\tau(\vartheta,\X_{T_{i^-}}),\\
         &\mathcal{L}_{t} ^{\text{move}}( \vartheta)=L(\vartheta,\X _{[T_{N_{t}},t]})\prod_{i=0}^{N_{t}-1} L(\vartheta,\X_{[T_{i},T_{i+1}[}).
      \end{align*}

    \end{theorem}
  
  \begin{rem}\label{remarque}
  The theoretical form of the likelihood allows, in principle, all characteristics of the process to depend on  a single parameter. In most applications however, the natural parametrisation involves distinct parameters for each characteristic, meaning that $\vartheta$ in Theorem~\ref{densite} typically takes the form 
  \begin{equation}\label{distinct}
  \vartheta=(\vartheta_{\beta},\vartheta_{k_\beta},\vartheta_{\delta},\vartheta_{k_\delta},\vartheta_{\tau},\vartheta_{k_\tau},\vartheta_{b}),
  \end{equation}
  where, for $\gamma\in\{\beta, \delta,\tau\}$,  the intensity $\gamma$ depends only on $\vartheta_\gamma$, the parameter $\vartheta_{k_\gamma}$ is exclusively involved in the transition kernel density $k_\gamma$, and $\vartheta_b$ appears only in the drift function of \eqref{EDSglob}. 
 In this setting, thanks to the product structure of the likelihood, each characteristic  can be estimated independently of the others. Corollary~\ref{corollaryMLE2} in the next section details the asymptotic behavior of the maximum likelihood estimator in this specific framework. 
   As an illustration,  in Section~\ref{sec:simus}, we  infer the birth kernel $k_\beta$ without specifying or estimating the remaining components of the model. 
  \end{rem}

\subsection{Local Asymptotic Normality}\label{subsec:LAN}
Intuitively, the LAN property of a parametric statistical model asserts that, in a neighbourhood of the true parameter $\vartheta$, the model becomes asymptotically equivalent to a Gaussian experiment. 
This fundamental notion, introduced in \cite{cam1960}, allows one  to reduce the statistical analysis of the model to a much simpler Gaussian framework. In particular, it entails both the asymptotic normality and the optimality of the maximum likelihood estimator; see for instance  \cite{ibragimov}. In this section, we show that under suitable regularity conditions, the BDMM process satisfies  the LAN property, and we clarify its consequences for the maximum likelihood estimator  in Corollaries~\ref{corollaryMLE} and~\ref{corollaryMLE2} .


We list below the  assumptions we make on the BDMM model  with parameter $\vartheta \in \Theta \subset \R^\ell$.  They are of two kinds: the first one \ref{hypH} ensures that the process is non-explosive and geometric ergodic. Conditions implying this hypothesis concern the intensity functions: it basically holds if the total intensity function $\alpha$ is bounded and the death intensity function $\delta$ sufficiently compensates  for the birth intensity $\beta$. Precise conditions are given in Appendix~\ref{appendix:ergodicity}, following \cite{emilien}. The other assumptions are regularity hypotheses on the way each feature of the BDMM process is parametrised. Although quite technical at first sight, we show in Section~\ref{lan examples} that they are satisfied in many situations, including the examples of Section~\ref{model examples}.

\begin{enumerate}[label=\textbf{(H)}, leftmargin=*]
\item The process is non-explosive, i.e. $N_t<\infty$ for any $t\geq 0$, the empty configuration $\text{\O}$ is a positive recurrent state,  and  there exist a measure $\mu_{ \infty}^\vartheta$ on $E$,  $0<r<1$ and $c : E \to (0,\infty)$, possibly depending on $\vartheta$,  such that for any $t>0$, any $\x_0 \in E$ and any measurable and bounded function $g $ on $E$,
$$ \displaystyle\left| \E_{\x_0}^\vartheta \left(g(\X_{t})\right)-\int_{E}g(\y)\mu_{\infty}^\vartheta(d\y)\right|\leq c(\x_0) \|g\|_{\infty} r^t.$$ \label{hypH}
\end{enumerate}
Concerning the intensity function $\gamma$ (whether $\gamma$ is $\beta$, $\delta$ or $\tau$), we assume that it is differentiable with respect to $\vartheta$ and the two following hypotheses. For a  function $f$, we denote by $\nabla f$ the gradient (or the Jacobian matrix when $f$ is vector-valued) with respect to $\vartheta$. When  $f$ depends on additional variables, $\nabla f$ will always refer to the gradient with respect to the first argument. We  also denote by $B(\vartheta,\rho)$ the ball centred at $\vartheta$ with radius $\rho>0$. 

\begin{enumerate}[label=\textbf{\Alph*($\gamma$)},leftmargin=*]
\item For all $\vartheta\in \Theta$,   $ \x\mapsto \frac{ \|\nabla \gamma( \vartheta,\x)\|^2}{ \gamma( \vartheta,\x)} \1_{\gamma(\vartheta,\x)>0}$ is bounded. \label{cond:A}

\item For all $ \vartheta\in \Theta$, $\exists \,\rho( \vartheta)>0$ such that  $ \x\mapsto f^ \gamma_{ \vartheta}(\x, \rho( \vartheta))\gamma( \vartheta,\x)$ is bounded, where \label{cond:B}
\end{enumerate}
\begin{align}\label{f_gamma}
f^{ \gamma}_{ \vartheta}(\x,\rho)=\displaystyle\sup_{u \in B( \vartheta, \rho) }\left(\left\| \nabla \gamma( u,\x)-\nabla \gamma( \vartheta,\x)\right\|^2\right)
\frac{\1_{\gamma(\vartheta,\x)>0}}{ (\gamma(\vartheta,\x))^2}.
\end{align}
Similarly, we assume that the transition kernel $k_\gamma$ is differentiable with respect to $\vartheta$ and the two following hypotheses. 

\begin{enumerate}[label={${\textbf{\Alph*}}$\textbf{(}$k_{ \gamma},\gamma$\textbf{)}},leftmargin=*]
\item For all $ \vartheta\in \Theta$,  $ \x\mapsto  \gamma( \vartheta,\x)  \int_{E}\frac{ \|\nabla k_{ \gamma}( \vartheta,\x,\y)\|^2}{k_{ \gamma}( \vartheta,\x,\y)}\1_{k_\gamma(\vartheta,\x,\y)>0}\nu_{ \gamma}(\x,d\y)$ is bounded. \label{cond:Atilde}
\item For all $\vartheta\in \Theta$,  $\exists\, \rho( \vartheta)>0$ such that $ \x\mapsto  f^{k_{ \gamma}}_{ \vartheta}(\x, \rho( \vartheta))\gamma( \vartheta,\x)$ is bounded, where
\label{cond:Btilde}
\end{enumerate}
\begin{align}\label{f_kgamma}
f^{ k_{ \gamma}}_{ \vartheta}(\x,\rho)=\int_{E} \sup_{u \in B( \vartheta, \rho) }\left( \left\| \nabla k_{ \gamma}( u,\x,\y)-\nabla k_{ \gamma}( \vartheta,\x,\y)\right\|^2  \right)\frac{\1_{k_\gamma(\vartheta,\x,\y)>0}}{k_{ \gamma}( \vartheta,\x,\y)} \nu_{ \gamma}(\x,d\y).
\end{align}
Finally, concerning the inter-jump motion given as a solution of \eqref{EDSglob}, we assume that its drift function $\bar b_n$ is differentiable with respect to $\vartheta$ and:
\begin{enumerate}[label=\textbf{\Alph*$_{\text{move}}$},leftmargin=*]
\item For all  $n\geq 1$,  $\forall \vartheta\in \Theta$,   $ \x\mapsto  \left\| \bar \sigma_n (\x)^{-1} \nabla \bar b_n( \vartheta,\x) \right\|$ is bounded. \label{cond:Amove}
\item   For all $n\geq 1$, $\forall \vartheta\in \Theta$,  $\exists\, \rho( \vartheta)>0$ such that $\x \mapsto \|\bar\sigma_n(\x)^{-1}\|^2 f^{\text{move}}_{ \vartheta}(\x, \rho( \vartheta))$  is bounded, where 
\label{cond:Bmove}
\end{enumerate}
\begin{align}\label{f_L}
f^{\text{move}}_{ \vartheta}(\x,\rho)= \displaystyle\sup_{u \in B (\vartheta, \rho) }\left\| \nabla \bar b_n( u,\x)-\nabla \bar b_n( \vartheta,\x) \right\|^2.
\end{align} 

In order to state the LAN property, we introduce,  for any $\vartheta\in\Theta$, the following $\mathcal F_{t}$-martingales, together with their associated angle brackets (see Appendix~\ref{appendix martingale} for background on these notions).  The martingale property and the explicit expressions of their brackets are established in the proof, detailed in  Appendix~\ref{proof LAN}, of the following theorem. We recall that $a_{n}= \bar\sigma_{n} \bar\sigma_{n}^T$. 
\begin{align}
M_{ \gamma}^{ \vartheta}(t)&= \frac{ 1}{ \sqrt{t}}\displaystyle \int_{0}^{t} \frac{ \nabla \gamma( \vartheta,\X_{s^-})}{ \gamma( \vartheta,\X_{s^-})} d\tilde N^{\gamma}_{s},\label{M_gamma}\\
M_{ k_\gamma}^{ \vartheta}(t)&= \frac{ 1}{ \sqrt{t}}\left(\displaystyle \int_{0}^{t} \frac{ \nabla k_\gamma ( \vartheta,\X_{s^-},\X_s)}{ k_\gamma ( \vartheta,\X_{s^-},\X_s)} dN^{\gamma}_{s}\right.\nonumber\\ 
&\left.\hspace{2cm}-\int_{0} ^{t} \int_{E} \nabla k_\gamma( \vartheta,\X_{s},\y)\nu_{ \gamma}(\X_{s},d\y) \gamma( \vartheta, \X_{s})ds\right),\label{M_kgamma}\\
M_{L}^{ \vartheta}(t)& = \frac{ 1}{ \sqrt{t}}\displaystyle \sum_{i\geq0}\int_{0}^{t} \1_{ [T_{i}; T_{i+1}[}(s) (\nabla \bar b_n( \vartheta,\X_{s^-}))^T a_{n}^{-1}( \X_{s^-}) d V^{(i)}_{s-T_i},\label{M_L}
\end{align}
where $\tilde N^{\gamma}_{t}=N^{\gamma}_{t}- \int_{0}^t \gamma(\vartheta, \X_{s})ds$ and  $V_s^{(i)}=\left(\Z^{(i)}_{s}-\Z^{(i)}_{0} -\int_0^{s}\bar b_n( \vartheta, \Z^{(i)}_{u})du\right)\1_{s\geq T_i}$ if we agree that  $\X_t=\{(\Z^{(i)}_{t-T_i},\m^{(i)})\}$ when $t\in[T_i,T_{i+1}[$. Their associated angle brackets are given by
\begin{align}
\langle M_{ \gamma}^{ \vartheta} \rangle(t)&
=\frac{ 1}{ t}\int_{0}^{t}\1_{\gamma>0}\frac{ (\nabla \gamma) (\nabla \gamma)^T}{ \gamma}( \vartheta,\X_{s})ds,\label{brackets_gamma}\\
\langle M_{ k_{ \gamma}}^{ \vartheta}\rangle(t)&=\frac{1}{t}\int_{0} ^{t} \int_{E} \1_{k_\gamma>0} \frac{ \nabla k_\gamma (\nabla k_\gamma)^T}{ k_\gamma}( \vartheta,\X_{s},\y)\nu_{ \gamma}(\X_{s},d\y) \gamma( \vartheta, \X_{s})ds,\label{brackets_kgamma}\\
\langle M_{ L}^{ \vartheta}\rangle(t)&=\frac{1}{ t}\displaystyle\int_{0}^{t}  (\nabla \bar b_n( \vartheta,\X_{s}))^T a_{n}^{-1}( \X_{s})\nabla \bar b_n( \vartheta,\X_{s}) ds.\label{brackets_L}
\end{align}
Finally denoting  
 $\mathcal S=\{L,\beta,\delta,\tau,k_\beta,k_\delta,k_\tau\}$, we introduce 
 \begin{align}\label{defM}
 M^{\vartheta}(t) =\sum_{i\in\mathcal S} M_i^{\vartheta}(t),
 \end{align}
 the angle brackets of which reads, as proved in Appendix~\ref{proof LAN}, 
 \begin{align}\label{defangleM}
 \langle M^{\vartheta}\rangle(t) = \sum_{i\in\mathcal S} \langle M_{ i}^{ \vartheta} \rangle(t). \end{align}

\begin{theorem}\label{LAN}
Under the same conditions as  in Theorem~\ref{densite}, 
suppose that \ref{hypH}, \ref{cond:A}, \ref{cond:B}, \ref{cond:Atilde}, \ref{cond:Btilde},  \ref{cond:Amove} and \ref{cond:Bmove} are satisfied (whether $\gamma$ is $\beta$, $\delta$ or $\tau$). 
Then for all $\vartheta \in \mathbb{R}^\ell$, $t>0$, and $h\in\R^\ell$, we have the following decomposition:
 \begin{align}\label{eqlambda}
\log \frac {\mathcal{L}_{t}(\vartheta+\frac{h}{\sqrt t})}  {\mathcal{L}_{t}(\vartheta)} =h^T M^{ \vartheta}(t)- \frac{ 1}{ 2}h^T\langle M^{ \vartheta} \rangle (t) h+Rem_{ \vartheta,h}(t),
 \end{align}
where $Rem_{ \vartheta,h}(t)\to 0$ in $\P_{\x_0}^\vartheta$-probability when $t \to+\infty$, for any $\x_0\in E$, and where  $M^{\vartheta}$ and $\langle M^{ \vartheta} \rangle$ are given by \eqref{defM} and \eqref{defangleM}. Moreover, for all $\vartheta\in\Theta$, 
\begin{align}\label{thmeq2}
\left( M^{ \vartheta}(t) ,\langle  M^{ \vartheta}\rangle (t)\right) \xrightarrow[t\to+\infty]{}(\mathcal{N}\left(0,tJ( \vartheta)),tJ( \vartheta)\right),
\end{align}
where the joint convergence takes place in distribution and where the matrix $J(\vartheta)$, of size $(\ell,\ell)$, has for expression:
\begin{align}
J( \vartheta)&=\lim_{t \to+\infty} \langle M^{ \vartheta}\rangle(t) \label{Jemp}\\
&= \int_E (\nabla \bar b_n( \vartheta,\x))^T a_{n}^{-1}( \x)\nabla \bar b_n( \vartheta,\x)  \mu_{\infty}^\vartheta(d\x) \nonumber\\
 &\hspace{0.5cm}+ \sum_{ \gamma\in \{ \beta, \delta, \tau\}} \int_E \1_{\gamma>0}\frac{ (\nabla \gamma) (\nabla \gamma)^T}{ \gamma}( \vartheta,\x) \mu_{\infty}^\vartheta(d\x)\nonumber \\
& \hspace{0.5cm}+ \sum_{ \gamma\in \{ \beta, \delta, \tau\}} \int_E \int_E   \1_{k_\gamma>0} \frac{ \nabla k_\gamma (\nabla k_\gamma)^T}{ k_\gamma}( \vartheta,\x,\y)\nu_{ \gamma}(\x,d\y) \gamma( \vartheta, \x) \mu_{\infty}^\vartheta(d\x). \nonumber
\end{align}
\end{theorem}

Under additional regularity assumptions on the model, as detailed for instance in \cite{ibragimov}, the LAN property implies  that the maximum likelihood estimator is asymptotically efficient and asymptotically Gaussian with covariance $J(\vartheta)^{-1}$. The following corollary summarise this consequence, the proof of which may be found in \cite{ibragimov}.
\begin{cor}\label{corollaryMLE}
In addition to the setting of Theorem~\ref{LAN}, assume that the regularity conditions N2-N4 of  \cite[Chapter III]{ibragimov} are satisfied. Then the maximum likelihood estimator $\hat\vartheta_t$, defined by $\L_t(\hat\vartheta_t)=\sup_{\vartheta} \L_t(\vartheta)$, where $\L_t$ is given in Theorem~\ref{densite},  is asymptotically efficient and satisfies 
$$\sqrt t (\hat\vartheta_t - \vartheta) \xrightarrow[t\to+\infty]{} \mathcal N(0, J(\vartheta)^{-1}),$$
where the convergence is in distribution and $J(\vartheta)$ is given by \eqref{Jemp}.
\end{cor}

In practice, the asymptotic covariance matrix $J(\vartheta)^{-1}$ is estimated by $J_t(\hat\vartheta_t)^{-1}$, where $\hat\vartheta_t$ is the maximum likelihood estimator and $J_t$ is the empirical version of $J$ given by \eqref{Jemp}, viz.
\begin{equation}\label{JNempir}
J_t(\vartheta)=\langle M^{\vartheta}\rangle(t),
\end{equation}
where $\langle M^{\vartheta}\rangle$ is given by \eqref{defangleM} and involves the terms \eqref{brackets_gamma}, \eqref{brackets_kgamma} and \eqref{brackets_L}. 

In connection with Remark~\ref{remarque}, which considers  the case where each characteristic of the process depends on distinct parameters, the following corollary summarises the empirical covariance and asymptotic distribution of the corresponding maximum likelihood estimators. 

\begin{cor}\label{corollaryMLE2}
Under the same assumptions as in Corollary~\ref{corollaryMLE}, assume that each characteristic of the process depends on distinct parameters, so that $\vartheta$ takes the form \eqref{distinct},  as described in Remark~\ref{remarque}. For $\gamma\in\{\beta, \delta,\tau\}$, consider the corresponding maximum likelihood estimators $\hat\vartheta_\gamma$, $\hat \vartheta_{k_\gamma}$ and $\hat \vartheta_b$. Then, as $t$ tends to infinity, 
$ \sqrt t \left(\langle M_{ \gamma}^{\hat\vartheta_\gamma} \rangle(t)\right)^{1/2} (\hat\vartheta_\gamma - \vartheta_\gamma)$, resp. 
$ \sqrt t \left(\langle M_{ k_{ \gamma}}^{ \hat\vartheta_\gamma,\hat\vartheta_{k_\gamma}}\rangle(t)\right)^{1/2} (\hat \vartheta_{k_\gamma}-\vartheta_{k_\gamma})$ and
$ \sqrt t \left(\langle M_{ L}^{\hat  \vartheta_b}\rangle(t)\right)^{1/2} (\hat \vartheta_b - \vartheta_b)$, converge in distribution to a standard Gaussian random variable.  Here, $\langle M_{ \gamma}^{\hat\vartheta_\gamma} \rangle(t)$ is given by \eqref{brackets_gamma} where the parameter $\vartheta$ appears only  through its component $\vartheta_{\gamma}$,  estimated by  $\hat\vartheta_{\gamma}$. Similarly, $\langle M_{ k_{ \gamma}}^{ \hat\vartheta_\gamma,\hat\vartheta_{k_\gamma}}\rangle(t)$ is given by \eqref{brackets_kgamma} and  involves only the components $\vartheta_\gamma$ and $\vartheta_{k_\gamma}$, while $\langle M_{ L}^{\hat  \vartheta_b}\rangle(t)$ is given by  \eqref{brackets_L} and depends solely on the parameter $\vartheta_b$. 
\end{cor}

 The results of this corollary are exploited in the numerical study and data analysis reported  in Section~\ref{sec:simus}, which focus on the estimation of $\vartheta_{k_\gamma}$, for $\gamma=\beta$, under a particular specification of the model.

\subsection{Examples}\label{lan examples}

We show in this section that the regularity conditions \ref{cond:A}, \ref{cond:B}, \ref{cond:Atilde}, \ref{cond:Btilde},  \ref{cond:Amove} and \ref{cond:Bmove} are satisfied under mild assumptions for the examples considered in Section~\ref{model examples}. 

\smallskip

For the intensity functions, whether $\gamma(\x)=\gamma$ is constant or $\gamma(\x)=\gamma n(\x)$, for $\gamma>0$, and the parameter of interest is $\vartheta=\gamma$, we have $f^{ \gamma}_{ \vartheta}=0$ in \eqref{f_gamma}, so that \ref{cond:B} is obviously satisfied. For \ref{cond:A}, it  clearly holds in the first case, while it is true in the second case if we assume that there is a maximal number of particles, i.e. $\sup_\x n(\x)<n^*$ for some $n^*>0$. Note that the latter restriction is purely theoretical, but it is not a limitation in practice. It is also an assumption made in  \cite{ronanfrederic} and it amounts to the specific setting of condition \eqref{eq30} presented in appendix,  that implies geometric ergodicity of the process, i.e. \ref{hypH}. 

\smallskip

Concerning the transition kernels, we focus below on the examples of birth and mutation kernels presented in Section~\ref{model examples}. For the death kernel, we considered only the uniform distribution as an illustrative example in that section. Although a parametrisation could be introduced for the death  kernel, we do not pursue this possibility in what follows.

For the birth kernel, remember that it reads as in \eqref{birth kernel}, i.e. $k_{ \beta}(\x, \y)=p_{m}k_{ \beta}^m(\x,z)$. Let us examine Examples~\ref{mixture} and \ref{potentiel} for   $k_\beta^m$.


\smallskip

\rex{mixture} In this example $k_\beta^m$ is a mixture of Gaussian distributions with deviation $\sigma>0$, so that the parameters of interest for the birth kernel are in this case  $p_m$, $m\in\M$, and $\sigma$. It is not difficult to check that if $\Lambda$ is a bounded set, then all partial derivatives of  $k_{ \beta}(\x,\y)$ with respect to each parameter are uniformly bounded in $\x$ and $\y$. On the other hand  $k_{ \beta}(\x,\y)$ itself is lower bounded in $\x$ and $\y$ under the same setting. Consequently, \ref{cond:Atilde} holds true whenever $\Lambda$ is a bounded set and the birth intensity $\beta(\x)$ is bounded in $\x$. The latter is in particular true if $\beta$ satisfies the setting discussed in the previous paragraph dedicated to the intensity functions.  For  \ref{cond:Btilde}, note that a rough control consists in upper-bounding the squared norm of the difference in \eqref{f_kgamma} by twice the sum of each squared norm. Then the same bounds as previously can be used to check  \ref{cond:Btilde} in the same setting. 

\smallskip

\rex{potentiel} Assume that $k_{ \beta}^m(\x,z) = e^{-\vartheta_m'S_m(\x,z)}/c^m(\x)$, 
where $\vartheta_m=(\vartheta_{m,m'})_{m'\in\M}$, $c^m(\x)=\int_\Lambda e^{-\vartheta_m'S_m(\x,z)}dz$ and $S_m$ is the vector 
$$S_m(\x,z)=\left(\sum_{\substack{x_i=(z_i,m_i)\in\x \\ m_i=m'}} \Phi_{m,m'}(z-z_i)\right)_{m'\in\M}.$$
This corresponds to a potential of the form $V(\x)=\sum_{i\neq j}\vartheta_{m_i,m_j} \Phi_{m_i,m_j}(z_i-z_j)$ in Example~\ref{potentiel}. The parameters of interest here are all the $\vartheta_{m,m'}$, $m,m'\in\M$. Assume that for all $m,m'\in\M$, $\Phi_{m,m'}$ is bounded by $L$. Then if we denote by $n_{m'}(\x)$ the number of particles in $\x$ having label $m'$, we have 
$$\|S_m(\x,z)\|^2\leq  \sum_{m'\in\M} L^2 n_{m'}^2(\x)  \leq  L^2 n(\x)  \sum_{m'\in\M} n_{m'}(\x) =  L^2 n^2(\x).$$ 
Assume in addition that there is a maximal number of particles $n^*$ and that $\Lambda$ is a bounded set, then using the previous upper-bound, we may prove by elementary inequalities that both \ref{cond:Atilde} and \ref{cond:Btilde} are satisfied. 

\medskip

\noindent Concerning the mutation kernel, let us  inspect Example~\ref{transitionmatrix}.

\smallskip

\rex{transitionmatrix} Remember that for this example we have
$$k_\tau(\vartheta,\x,\y) = \frac 1 {n(\x)} p_{m_i,m}$$
if there exists $x_i\in\x$ and $m\in \M$ such that $\y=(\x\backslash (z_{i},m_i))\cup (z_{i},m)$. 
The parameter of interest $\vartheta$  in this example corresponds to  the non-null entries $p_{m,m'}$ of the transition matrix. Denote by $p_*$ the minimal value of these non-null entries. We have $\|\nabla k_{ \tau}(\vartheta,\x,\y)\|^2=1/n(\x)^2$, so that  
\begin{align*}
\int_{E} \frac{ \|\nabla k_{ \tau}( \vartheta,\x,\y)\|^2}{k_{ \tau}( \vartheta,\x,\y)}& \1_{k_\tau(\vartheta,\x,\y)>0}\nu_{ \tau}(\x,d\y) \\
& =\int_{E}\frac{ \|\nabla k_{ \tau}( \vartheta,\x,\y)\|^2}{k_{ \tau}( \vartheta,\x,\y)^2}k_{ \tau}( \vartheta,\x,\y)\1_{k_\tau(\vartheta,\x,\y)>0} \nu_{ \tau}(\x,d\y) \\
&\leq \frac 1 {p_*^2} \int_{E} k_{ \tau}( \vartheta,\x,\y) \nu_{ \tau}(\x,d\y) = \frac 1 {p_*^2}.
\end{align*}
Consequently \ref{cond:Atilde} (where $\gamma=\tau$) is satisfied in this setting whenever $\tau(\x)$ is bounded in $\x$, which holds true in the setting discussed in the first paragraph of this section. On the other hand, since $\nabla k_{ \tau}(\vartheta,\x,\y)$ does not depend on $\vartheta$, $f^{ k_{ \tau}}_{ \vartheta}(\x,\rho)=0$ in  \eqref{f_kgamma}, so that \ref{cond:Btilde} trivially holds.

\medskip

Concerning the inter-jump motion and assumptions  \ref{cond:Amove} and \ref{cond:Bmove}, we consider as an example the Langevin diffusion introduced in Example~\ref{langevin}.

\smallskip

\rex{langevin} Assume the same kind of parametrisation of the potential as in Example~\ref{potentiel} above, that is $b_{i,n}(\Z,\m)=- \sum_{j\neq i}  \vartheta_{m_i,m_j} \Phi_{m_i,m_j}(z_{i}-z_{j})$ in \eqref{EDSglob},  where the parameters of interest are again all the $\vartheta_{m,m'}$, $m,m'\in\M$. Remember that $\sigma_i(\Z,\m)= \sigma_{m_{i}}$ is known. Since the drift function is linear in the parameters, it is easily seen that  \ref{cond:Amove} is satisfied if $\Phi_{m,m'}$ is  bounded for any $m,m'\in\M$ and there is a maximal number of particles $n^*$. As to  \ref{cond:Bmove}, we clearly have $f^{\text{move}}_{ \vartheta}(\x,\rho)=0$ in \eqref{f_L}, so that this condition holds trivially true.

\section{Simulations and data analysis}\label{sec:simus}

\subsection{Simulation study} \label{simu}
We simulate in this section several trajectories of the BDMM process, for characteristics explained below, and we evaluate the quality of estimation of a parameter involved in the birth kernel. The chosen characteristics are motivated by the application that we will conduct in the next section, and are in line with previous studies by \cite{briane} and \cite{ronanfrederic}. We present them below without precisely detailing each value.

The particles are observed on the square $\Lambda=[-10,10]^2$. They do not possess continuous marks, but can be of 6 different types, which correspond on one hand to 2 distinct particle natures (playing the role of Rab-11 type proteins and Langerin type proteins in the subsequent application), and on the other hand to 3 possible motion regimes: Brownian, sub-diffusive (according to an Ornstein-Uhlenbeck process), or super-diffusive (according to a Brownian motion with linear drift). Thus, the space $\M$ contains 6 different labels, but mutations will only be possible between motion regimes and not between the type of particles.

\begin{figure}[h]
\begin{center}
       \includegraphics[width=7.5cm]{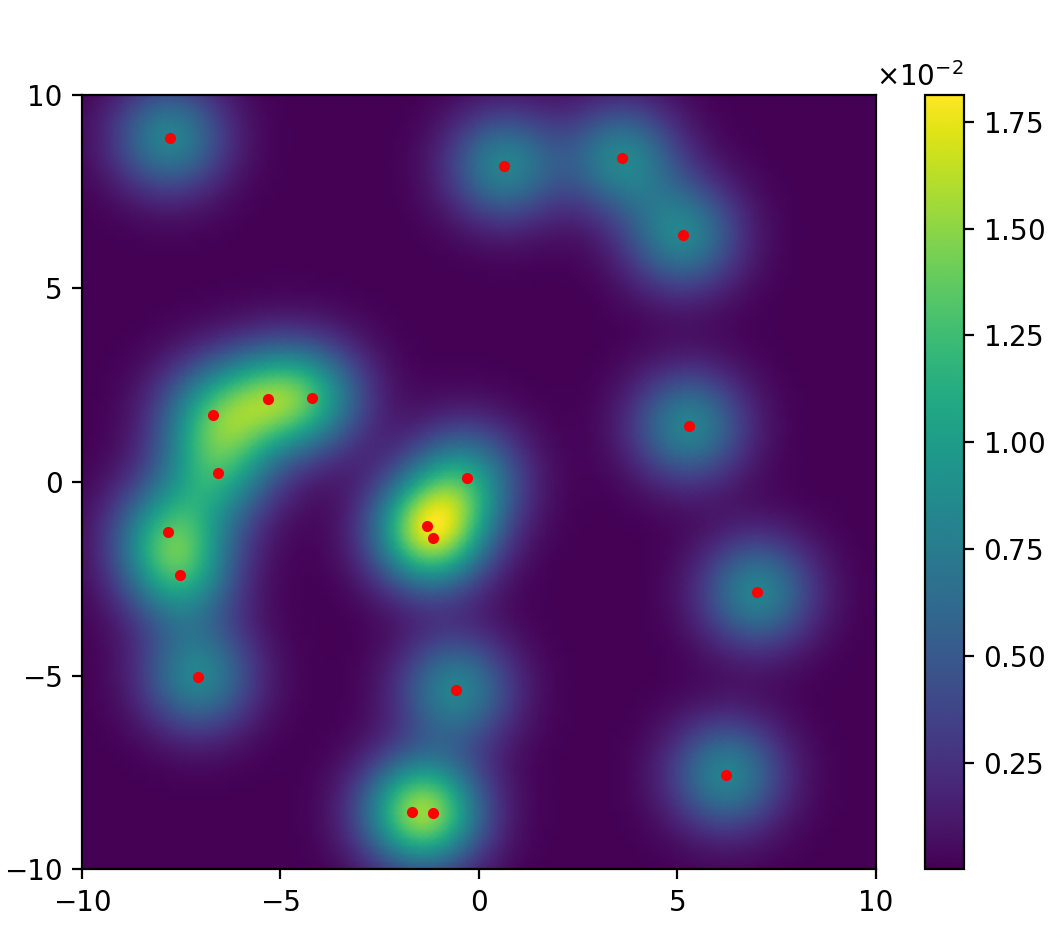} \\
\end{center}
        \caption{Example of the density \eqref{nucleus_naissances} for the birth of a new Langerin particle, given the Rab-11 particles already present (shown as red dots).}
        \label{figdensity}
    \end{figure}

We assume that the movement of each particle is independent of the others and corresponds to one of the three aforementioned regimes. 
 Birth rates are constant and differ only according to the Langerin or Rab-11 type. Death rates (also distinct for Langerin and Rab-11) are proportional to the number of particles present. Mutations correspond to a change in movement regime: their intensity is also constant. Regarding transitions, we assume that the death kernel is uniform: each particle has the same probability of disappearing when a death event occurs. For mutations, we follow the setting of Example~\ref{transitionmatrix}, i.e., these occur for any particle uniformly, according to a fixed transition matrix. Finally, for births, we assume that they generate an equiprobable Langerin or Rab-11 particle, with an equiprobable motion regime. The Rab-11 type particles are then generated uniformly in $\Lambda$. 
For the Langerin type particles, we assume that given the pre-jump configuration $\x\in E$, they are generated according to the density $k_\beta^L(\x,.)$ defined for all $z\in\R^2$ by
\begin{equation}\label{nucleus_naissances}
k_{ \beta}^{L}(\x,z)= \frac{ p}{ n_R(\x)}\sum_{i=1}^{n_R(\x)} \frac{ 1}{ 2\pi\sigma^2}\exp\left( \frac{-\|z-z^R_i \|^2}{ 2\sigma^2}\right)+\frac{ (1-p)}{ | \Lambda|}\1_{ \Lambda}(z),
\end{equation}
where $n_R(\x)$ denotes the number of existing Rab-11 particles and $z_i^R$ their spatial coordinates. 
This density is parametrized by $p\in [0,1]$ and $\sigma>0$. It is a convex combination between a mixture of normal distributions centered around the existing Rab-11 particles, with variance $ \sigma^2$, like in Example~\ref{mixture}, and a uniform distribution over $\Lambda$. The closer the parameter $p$ is to $1$, the more likely Langerin and Rab-11 particles will tend to be colocalized, and the smaller the standard deviation $ \sigma$ is, the stronger the proximity between the two types of particles will be. An illustration of this density is shown in Figure~\ref{figdensity}, for a given configuration of Rab-11 particles represented by red dots.

Note that the above birth density is defined on $\R^2$ and not only on $\Lambda$, meaning that some Rab-11 particles might give birth to Langerin particles outside $\Lambda$. Also, even if all Rab-11 particles are generated inside $\Lambda$, some may exit $\Lambda$ during their motion, and may further give birth to Langerin particles inside or outside $\Lambda$. This means that the process we are generated is actually defined on $\R^2$, even if most of the particles lie on $\Lambda$. But for inference, we assume that we observe the process on $\Lambda$ only. The situation is therefore slightly different from the (ideal) theoretical setting of Theorem~\ref{densite} where the likelihood is given for a process both defined and observed on $\Lambda$.
We chose this non-ideal setting to be in line with the application of the next section, where the observed particles are near the membrane of a living cell, but might exit the focal depth of the microscope during their motion (thus disappear), and hidden particles might ``give birth'' to observed particles. 
However, we expect that these border effects are negligible for the births in $\Lambda$, whenever $p$ and $\sigma$ are small, a situation where the interactions between the particles are rare and short range, as expected in the subsequent application. 

We are interested in the estimation of the parameter $(p,\exp \sigma)$ involved in the birth transition density $k_\beta^L$, from the observation of a single trajectory of the process on $\Lambda$ and during the  time interval $[0,t]$. 
 As emphasized in Remark~\ref{remarque}, we do not need to know all features of the process to carried out this estimation and we simply estimate $(p,\exp \sigma)$ by 
\begin{equation}\label{thetahat}
(\hat p,\widehat{\exp \sigma})=\underset{(p,\exp \sigma)}{\rm{argmax}} \sum_{i\in \mathcal{B}^L_{t}}\log k_{ \beta}^L(\X_{T_{i^-}^L},z_{T_{i}^L}),
\end{equation}
where $\mathcal{B}^L_{t}$ denotes the set of birth times for the Langerin particles observed on $\Lambda$ up to time $t$,  $(T_i^L)_{i\in \mathcal{B}^L_{t}}$ are these birth times, and $(z_{T_{i}^L})_{i\in \mathcal{B}^L_{t}}$ are the coordinates of the new Langerin particles in $\Lambda$.

\begin{figure}
\begin{center}
             \includegraphics[height=3.9cm]{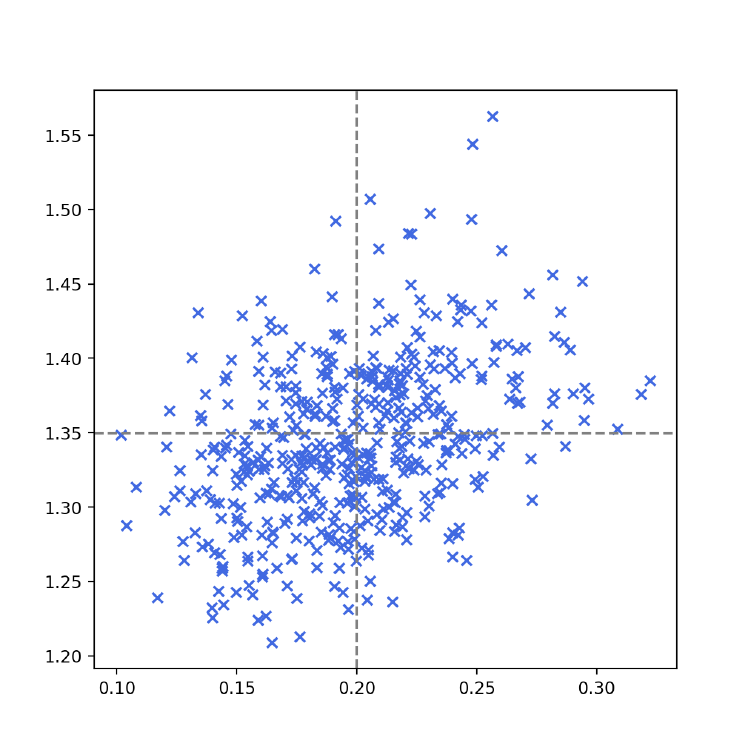}
             \includegraphics[height=3.9cm]{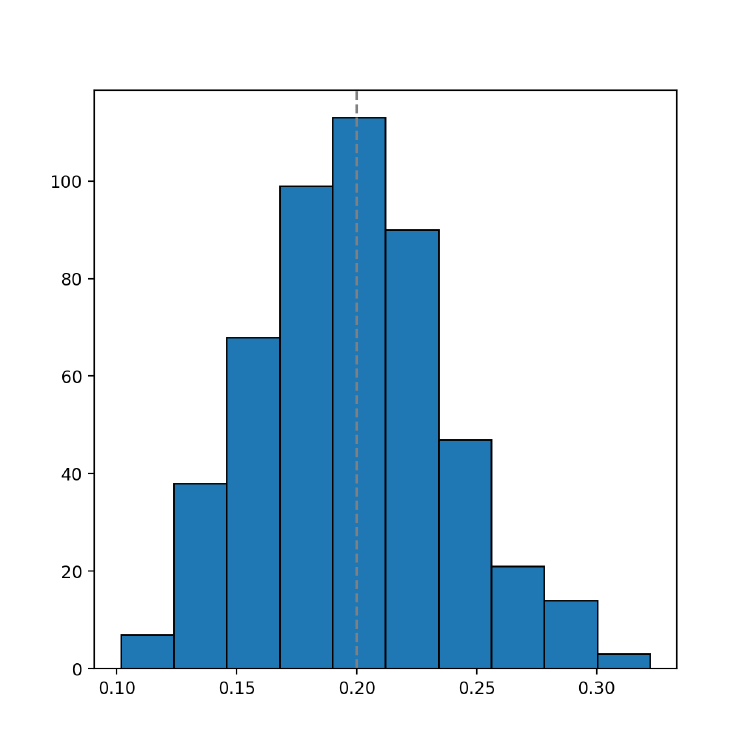}
             \includegraphics[height=3.9cm]{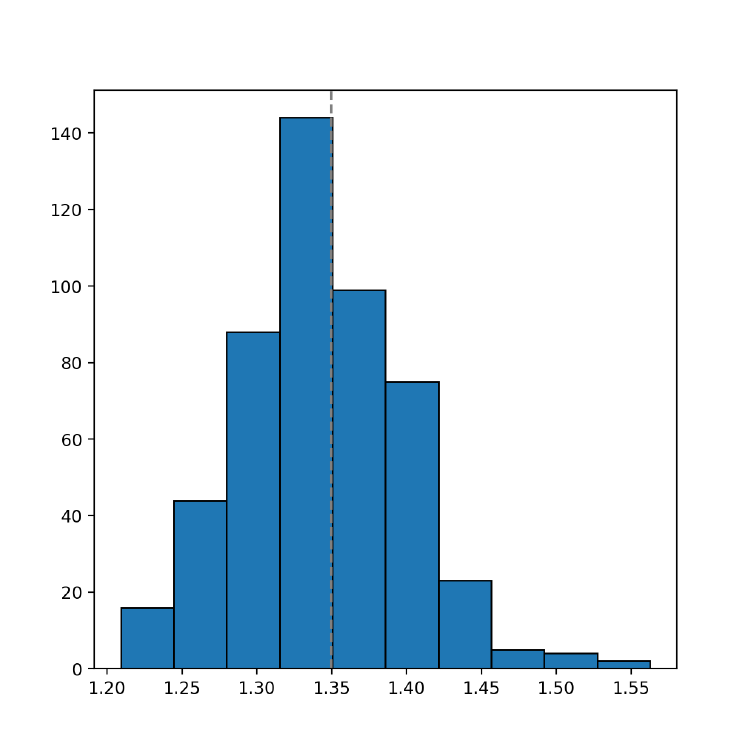}
\end{center}
        \caption{Estimation results for $(p,\exp \sigma)=(0.2,1.35)$ by MLE, based on 500 simulations of the BDMM process specified in Section~\ref{simu}. Left: scatterplot of all 500 estimates $(\hat p, \widehat{\exp \sigma})$. Middle: marginal distribution of $\hat p$. Right: marginal distribution of $\widehat{\exp \sigma}$.}
        \label{figellipse}
    \end{figure}

We performed 500 simulations of trajectories for $(p,\exp \sigma)=(0.2,1.35)$, on a time interval similar to the data studied in Section~\ref{section:simus:appli}. 
Figure~\ref{figellipse} displays the full set of 500 estimates obtained from \eqref{thetahat}, together with their marginal distributions. These empirical results appear consistent with the asymptotic Gaussian law described in Corollary~\ref{corollaryMLE2}.
In each case, we also computed the asymptotic covariance matrix 
$(\langle M_{ k_{ \beta}}\rangle(t))^{-1}$ of Corollary~\ref{corollaryMLE2}, given by \eqref{brackets_kgamma}. This matrix depends on the unknown parameter $(p,\exp \sigma)$,  which is replaced by its estimate, as well as on the birth intensity for the Langerin particles. In this example, the latter is assumed to be constant and is therefore simply estimated by the number of Langerin births observed over $[0,t]$ divided by $t$. 
 The corresponding 95\% confidence intervals yielded  empirical coverages of $95.0\%$ for $p$, $94.6\%$ for ${\exp \sigma}$, and $96.0\%$ for the joint confidence ellipsoids of both parameters.

Concerning computational complexity, each trajectory involved approximately 200 birth events, with about 20 to 30 Rab-11 particles present at each birth. The estimation of $\vartheta$ then took about 1 second per trajectory on a standard laptop, using standard optimization functions available in Python. In view of \eqref{thetahat}, this computation cost is linear in the number of observed births and in the number of Rab-11 particles.  Note that this optimization posed no particular problems, since the objective function in \eqref{thetahat} is bell-shaped and present no local maxima (see an instance in Figure~\ref{figvrs}). As a result, taking different initial values for the  optimization routine, even outside the ranges showed in Figure~\ref{figellipse}, gave the same results. The computation of the covariance matrix  is in turn more involved, since it requires to evaluate the integral \eqref{brackets_kgamma}, and took about 5 seconds per trajectory.

\subsection{Data analysis}\label{section:simus:appli}

\begin{figure}
  \centering
     \centering \includegraphics[height=5cm]{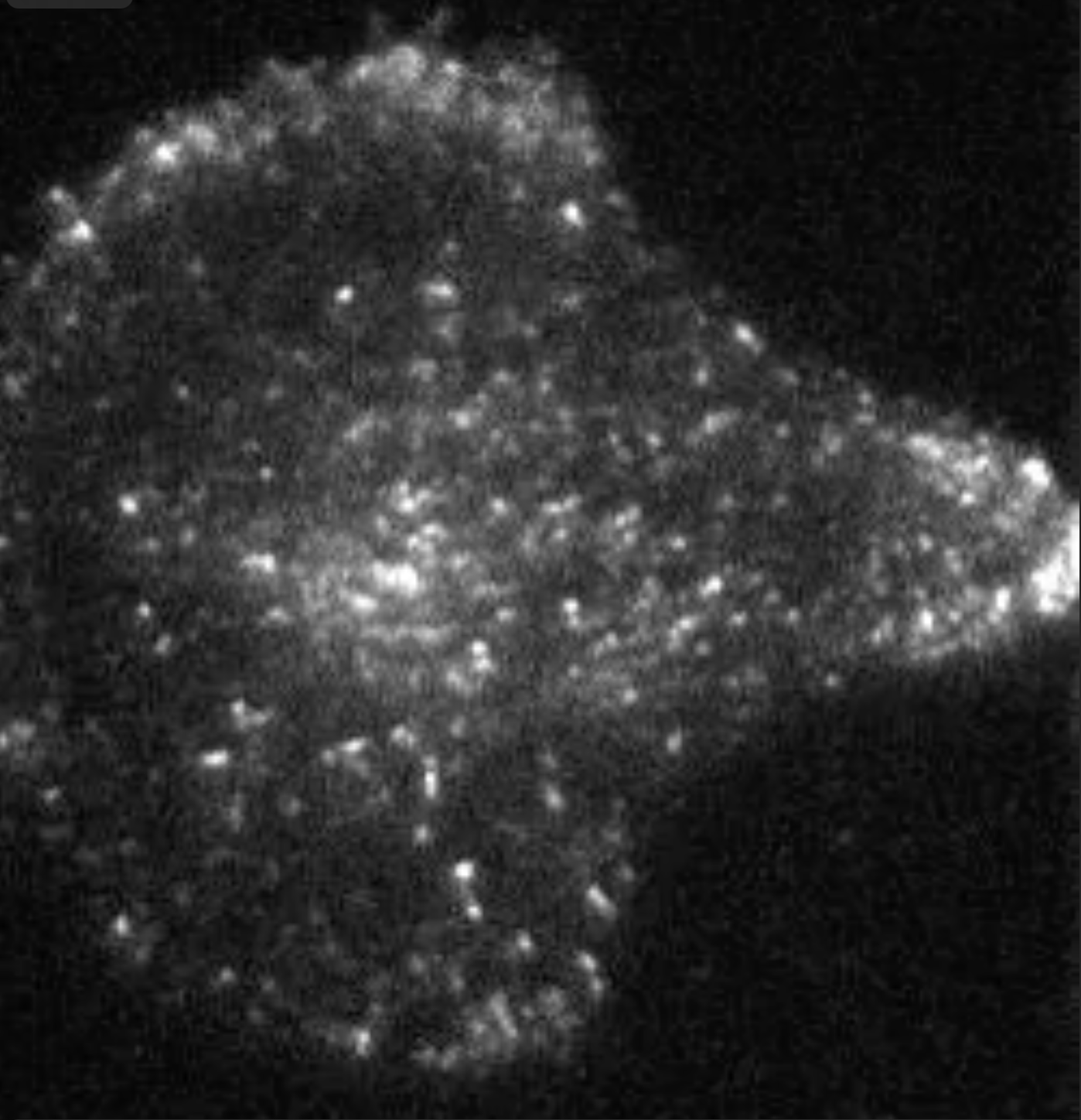}
  \hspace{50pt}
     \centering \includegraphics[height=5cm]{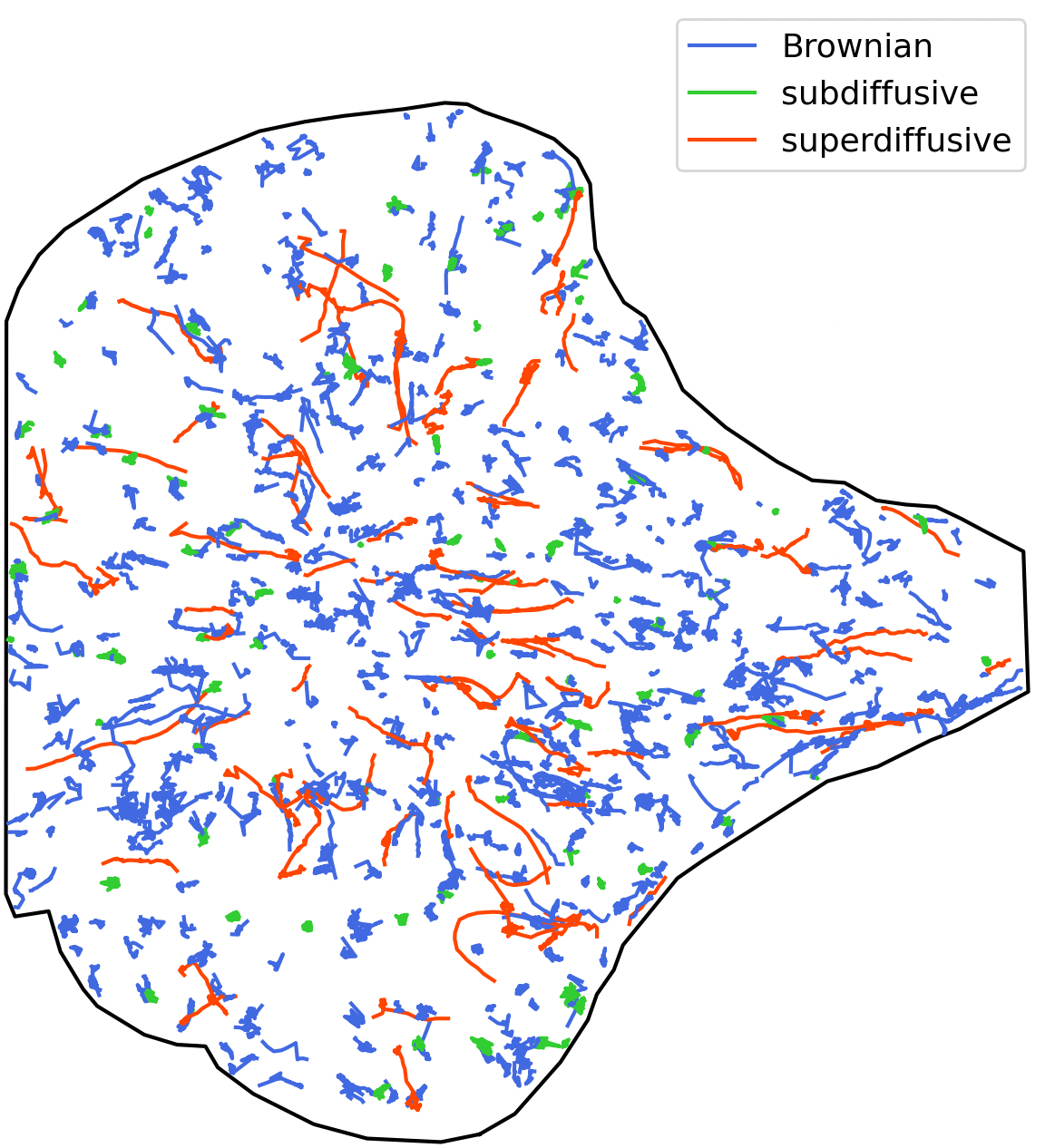}
  \caption{Left: a frame from the raw video sequence analysed in Section~\ref{section:simus:appli}, showing Langerin proteins. Right: superposition of all trajectories of Langerin proteins obtained after post-processing (i.e. segmentation and tracking) of the whole raw video. The colors indicate the estimated motion regime.}
   \label{real_data}
\end{figure}

The dataset we consider comes from the observation by TIRF (Total Internal Reflection Fluorescence) microscopy of the intracellular traffic of some molecules near the membrane of a living cell \cite{Boulanger2014}. This provides a video sequence showing two types of proteins observed simultaneously in the same cell: Langerin proteins and Rab-11 proteins. The image on the left of Figure~\ref{real_data} shows a frame obtained from the video of the Langerin proteins. After post-processing following  \cite{pecot2014background}, the proteins of interest are identified, represented by a point and tracked by the U-track algorithm  \cite{jaqaman2008} along the video sequence to provide trajectories, such as those visible in the right representation of Figure~\ref{real_data}. These trajectories have been further analysed by the method developed in \cite{Briane18} to classify them into three diffusion regimes  (these are the colors visible on the figure). 
The full  dynamics is in line with a BDMM process: in addition to their displacement, some proteins disappear in the course of time, others appear, and finally some change their diffusion regime, which corresponds in our model to a mutation (even if it is not a mutation in the biological sense of the term).

\begin{figure}
        \centering \includegraphics[scale=0.58]{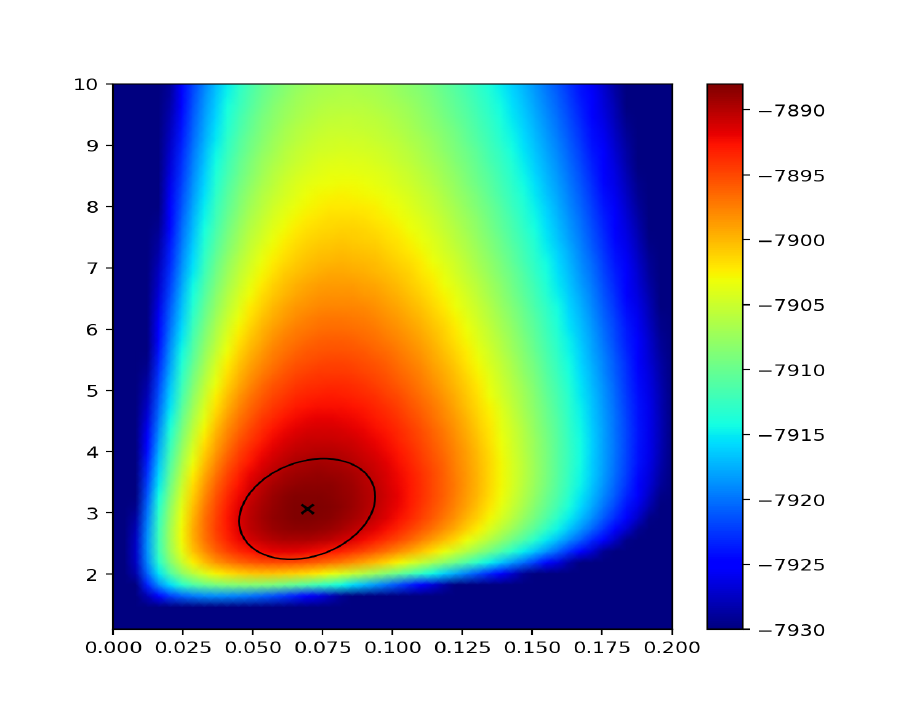}
        \caption{Likelihood of the BDMM model for the considered dataset, as a function of the parameters $p$ and $\exp \sigma$. The estimated $95\%$ confidence ellipsoid, centred at the maximum likelihood, is overlaid in black. }
        \label{figvrs}
    \end{figure}

An important biological question is to determine if Langerin proteins tend to appear close to existing Rab-11 proteins. This is the phenomenon of colocalization, indicating a strong interaction between these two types of proteins. 
To answer this question we estimate by maximum likelihood the birth kernel of the Langerin proteins, assuming that it is of the form \eqref{nucleus_naissances}, as in our simulation study. We obtained the estimate $(\widehat p, \widehat{\exp \sigma})=(0.07,3)$, where the unit of $\sigma$ is the pixel (the image being of size $250 \times 283$, each pixel representing an area of $80\times 80$ nm$^2$). Figure~\ref{figvrs} shows the likelihood value for our dataset, with respect to $p$ and $\exp\sigma$, along with the empirical  $95\%$ confidence ellipsoid around the maximum, obtained from the estimation of the covariance matrix, see  Corollary~\ref{corollaryMLE2}. This estimate suggests that, for this dataset, about $7\%$ of Langerin proteins are colocalized with Rab-11 ones, this proportion being significantly positive in view of the confidence ellipsoid, and  each colocalized Langerin protein appears (with high probability) within $3\widehat\sigma\times 80=264$ nm of a Rab-11 protein, where $\widehat\sigma=\log(\widehat{\exp \sigma})$.

\section{Discussion}\label{discussion}

Our contribution addressed the parametric inference of the birth-death-move (BDM) process with possible mutations. This process involves several characteristics: i) the motions of particles; ii)  the intensities of births, deaths, and  mutations; iii) the transition kernels for the births, for the deaths and for the mutations. We  derived  the full likelihood of the process in  Theorem~\ref{densite}, and provided theoretical guarantees for maximum likelihood estimation, including its efficiency and asymptotic normality with explicit variance  (see Corollaries~\ref{corollaryMLE} and~\ref{corollaryMLE2}). In theory, the parameter can relate to any or even all characteristics of the process, see Remark~\ref{remarque}. 
However, estimating all of these characteristics in practice seems ambitious. For example, the sole parametric inference of the particles' motions amounts to estimating a system of stochastic differential equations, which is  itself the subject of intensive work in the literature. For this task, it is well known that the quality of estimation depends on the length of the observed trajectories, the chosen parametrisation and also the signal-to-noise ratio. We, of course, face the same limitations. Estimating the other characteristics of the BDM process may be easier, as illustrated in our application in Section~\ref{sec:simus}, where we estimated the birth kernel using a simple parametrisation. Fortunately, the product form of the likelihood allows us to estimate one characteristic without requiring knowledge or estimation of the others. In this setting, our theoretical results enable us to provide a confidence region for the parameter of interest, which properly accounts for the entire generative process, even if we are only interested in a specific characteristic.

Nevertheless, the statistical study of the BDM process is  far from exhaustive. Nonparametric estimation of its characteristics remains an important step, as it is often the first stage towards specifying a parametric model. This issue has been studied in \cite{ronanfrederic} for estimating the intensity functions, but it remains an open question for the estimation of the transition kernels. Another perspective involves the study of discrete-time observations, a more realistic framework than the ideal continuous-time setting considered in this contribution. Additionally, a limitation of our approach is that, in the likelihood of Theorem~\ref{densite}, all jump times  $T_i$ are assumed to be observed. While this assumption is generally not restrictive for the birth times and the death times (since birth or death of a particle is easily visible on a video), it can be unrealistic for mutation times. For example, in our application in Section~\ref{sec:simus}, a mutation represents a change in the diffusion regime of a particle, which is a typical difficult transition to observe. Since we were only interested in the birth kernel, these transition times did not need to be observed. But to address the full estimation problem, including the estimation of the mutation times, hidden Markov models must be considered. This is the subject of ongoing research.

\appendix

\section{Proof of Theorem~\ref{densite}}\label{appendix:likelihood}
To simplify notation in the proof, we drop the dependence in $\vartheta$ of all terms. 
To get the likelihood, we view $(\X_{t})_{t \geq 0}$ as a jump-move process, a more 
general Markov process than the BDMM process. A jump-move process alternates motions and jumps in $E$, see \cite{emilien}. It is defined through an intensity function $\alpha$, a jump kernel $K(\x,.)$ and a continuous Markov process $\Y$ on $E$ that drives the motion of the process. For a BDMM, we have $ \alpha(\x)= \beta(\x)+ \delta(\x)+ \tau(\x)$ and 
\begin{equation}\label{defK}
K(\x,.)=\frac{\beta(\x)}{\alpha(\x)}K_\beta(\x,.)+\frac{\delta(\x)}{\alpha(\x)}K_\delta(\x,.)+ \frac{  \tau(\x)}{ \alpha(\x)}K_{ \tau}(\x,.).\end{equation}
By our assumptions 
\begin{equation}\label{K_JM} K( \x,d\y)=k(\x,\y)\nu(\x,d\y),\end{equation}
where for any $F\in\mathcal E$, $\nu(\x,F)=\nu_\beta(\x,F)$ if there exist $A\subset \Lambda$ and $ \mathcal{I}\subset  \M$ such that $F=\x\cup (A \times \mathcal{I})$,  $\nu(\x,F)=\nu_\delta(\x,F)$ if $\x\setminus x_i\in F$ for some $i$, and   $\nu(\x,F)=\nu_\tau(\x,F)$ if there exist $m\in \M$ and $i$ such that $(\x\backslash (z_{i},m_i))\cup (z_{i},m)\in F$. 
In turn, $k(\x,\y)=k_\beta(\x,\y)\beta(\x)/\alpha(\x)$ if there exists $(z,m)$ such that $\y=\x\cup(z,m)$, $k(\x,\y)= k_\delta(\x,\y)\delta(\x)/\alpha(\x)$ if $\y=\x\setminus x_i$ for some $i$, and $k(\x,\y)= k_\tau(\x,\y)\tau(\x)/\alpha(\x)$ if there  exist $m\in \M$ and $i$ such that $\y=(\x\backslash (z_{i},m_i))\cup (z_{i},m)$.

Note that a jump-move process observed continuously on $[0,t]$ has by definition a cadlag trajectory on $E$ and can equivalently be described by the vector $\Xx$ defined by: 
 \begin{align}\label{Xx}
   \Xx=\displaystyle\sum_{n\in \N}(\X_0,\X_{]0,T_{1}[},T_{1},\X_{T_{1}},\X_{]T_{1},T_{2}[},\dots,T_{n}-T_{n-1},\X_{T_{n}},\X_{]T_{n},t[})  \1_{\left\{ {N_t=n}\right\}}.
    \end{align}
We set  for any $i \geq 1$, $S_{i}=T_{i}-T_{i-1},$ $\V_{i}=\X_{T_{i}},$ and we recall that 
$(\X_{t})_{T_{i}\leq t <T_{i+1}}=(\Y^{(i)}_{t})_{0\leq t\leq S_{i+1}}$, where $\Y^{(i)}$ is a continuous Markov process on $E$, identically distributed as $\Y$. 
With these notations $\Xx$  is equivalent to $\tilde\Xx$ with:
\begin{align}\label{Xx2}\tilde\Xx=\displaystyle\sum_{n\in \N}(\V_0,\Y^{(0)}_{[0,S_{1}]},S_1,\V_1,\Y^{(1)}_{[0,S_{2}]},\dots, S_n,\V_n,\Y^{(n)}_{[0,t-\sum_{j=1}^nS_j]})\1_{\left\{ {N_t=n}\right\}}.\end{align}
Note that with this formalism, $\Y^{(i)}_0=\V_i$ for all $i$, and $\tilde\Xx$  takes values in the space $\bigcup_{n\geq 0}{\mathcal W}_n$ where 
$${\mathcal W}_n= \bigotimes_{i=0}^n (E\times \mathcal C_{[S_i]}\times (0,\infty))  \times  E\times \mathcal C_{[t-\sum_{i=0}^n S_i]},$$
and where  $\mathcal{C}_{[s]}$ denotes the space of continuous functions from  $[0,s]$ to $E$.

As detailed in Section~\ref{exmove}, the continuous Markov process $\Y^{|n}_s=\Pi_n((\Z_s,\m))$ defined on $E_n$ and observed on $[0,t]$, given that $\Z_0=\z$, admits the density $L(\Z_{[0,t]},\m)$ with respect to the reference process $\Pi_n((\U_t,\m))$ where $\U_t=\z+\int_0^t \bar\sigma_n(\U_s,\m) dB_s$. 
Let us denote by $P^\z_{[t]}$ the distribution of this reference process observed on $[0,t]$. For any measurable and bounded function $h$ defined on $E_n$, we thus have 
 $\E(h(\Y^{|n}_{[0,t]}))=\int_{\mathcal C_{[t]}} h(\y) L(\y_{[t]}) dP^{\Y^{|n}_0}_{[t]}(\y),$ where $L(\y_{[t]})$ is a short notation for $L(\z_{[0,t]},\m)$ and where $\y_t=(\z_t,\m)$ for any $t\geq 0$. 
 Similarly, for $(\Y_t)_{t\geq 0}$ on $E$ defined by $\Y_t=\sum_{n\geq 0} \Y^{|n}_t \1_{\{\Y_0\in E_n\}}$, we will write without ambiguity 
 $$\E(h(\Y_{[0,t]}))=\int_{\mathcal C_{[t]}} h(\y) L(\y_{[t]}) dP^{\Y_0}_{[t]}(\y),$$ where the cardinality $n$ defining $P^{\Y_0}_{[t]}$  is implicitly given by $n=n(\Y_0)$.

Let $h$ be a measurable and bounded function. For $n\in \N$, we use the conditional law of $S_{n+1}$ given $(\X_t)_{t\leq T_{n}}$, $T_{n}$ and $\Y^{(n)}$, see \eqref{jump law}, to compute:
\begin{align*}
    \E  [&h(\tilde\Xx)\1_{\left\{ {N_t=n}\right\}} ]\\
    &=\E\left[h(\tilde\Xx)\1_{\left\{ {S_{n+1}>t-\sum_{j=1}^nS_j}\right\}}\1_{\left\{ {\sum_{j=1}^nS_j\leq t}\right\}}\right] \\
     &=\E\left[\E\left[h(\tilde\Xx)\1_{\left\{ {S_{n+1}>t- \sum_{j=1}^nS_j}\right\}}\1_{\left\{ {\sum_{j=1}^nS_j\leq t}\right\}}\big |(\X_t)_{t\leq T_n}, {T_n},\Y^{(n)}\right]\right]\\
    &=\E\left[h(\tilde\Xx)\1_{\left\{ {\sum_{j=1}^nS_j\leq t}\right\}}\exp\left(-\displaystyle\int_0^{t-\sum_{j=1}^nS_j}\alpha(\Y^{(n)}_u)d u\right)\right] \\
     &=\E\left[\E\left[h(\tilde\Xx)\1_{\left\{ {\sum_{j=1}^nS_j\leq t}\right\}}\exp\left(-\displaystyle\int_0^{t-\sum_{j=1}^nS_j}\alpha(\Y^{(n)}_u)d u\right)\big | (\X_t)_{t\leq T_n},{T_n}\right]\right].
\end{align*}
Since given $S_1,\dots,S_n$ and  $\Y^{(n)}_0=\V_n$, the Markov process $\Y^{(n)}$ observed on the time interval $[0, t-\sum_{j=1}^nS_j]$ 
admits the density $L(\y_{[t-\sum_{j=1}^nS_j]})$ with respect to $P^{\V_n}_{[t-\sum_{j=1}^nS_j]}$, we have 
\begin{align*}
   \E  [h(\tilde\Xx)&\1_{\left\{ {N_t=n}\right\}} ]\\
    =\E& \left[\displaystyle\int_{\mathcal C_{[t-\sum_{i=1}^n S_j]}} h(\V_0,\Y^{(0)}_{[0,S_{1}]},\dots, S_n, \V_n,\y^{(n)}_{[t-\sum_{i=1}^n S_j]})\1_{\left\{ {\sum_{j=1}^nS_j\leq t}\right\}}\right.\\
&\left.\times\exp\left(-\displaystyle\int_0^{t-\sum_{j=1}^nS_j}\alpha(\y^{(n)}_u)d u\right) L\left(\y^{(n)}_{[t-\sum_{i=1}^n S_j]}\right)dP^{\V_n}_{[t-\sum_{i=1}^n S_j]}(\y^{(n)})\right].
\end{align*}
Now, if we condition on $(\X_t)_{t\leq T_{n-1}}$, $T_{n}$ and $\Y^{(n-1)}$ we know that the post-jump location $\V_n$ admits the density $k(\X_{T_n^-},\vv)=k(\Y^{(n-1)}_{S_n},\vv)$ with respect to the measure $\nu(\X_{T_n^-},.)=\nu(\Y^{(n-1)}_{S_n},.)$ on $E$. Hence  
\begin{align*}
   \E  [h(\tilde\Xx)&\1_{\left\{ {N_t=n}\right\}} ]\\
  =\E&\bigg[\displaystyle\int_E \int_{\mathcal C_{[t-\sum_{i=1}^n S_j]}} h(\V_0,\Y^{(0)}_{[0,S_{1}[},\dots, S_n, \vv_n,\y^{(n)}_{[t-\sum_{i=1}^n S_j]})\1_{\left\{ {\sum_{j=1}^nS_j\leq t}\right\}}\\
  &\times \exp\left(-\displaystyle\int_0^{t-\sum_{j=1}^nS_j}\alpha(\y^{(n)}_u)d u\right)
L\left(\y^{(n)}_{[t-\sum_{i=1}^n S_j]}\right)dP^{\vv_n}_{[t-\sum_{i=1}^n S_j]}(\y^{(n)}) \\
&\times k(\Y^{(n-1)}_{S_n},\vv_n)\nu(\Y^{(n-1)}_{S_n},d\vv_n)\bigg].
\end{align*}
Conditioning on $(\X_t)_{t\leq T_{n-1}}$, $T_{n-1}$ and $\Y^{(n-1)}$, we use again the law of $S_n$ given by \eqref{jump law} to get 
\begin{align*}
   \E&  [h(\tilde\Xx)\1_{\left\{ {N_t=n}\right\}} ]\\
   =&\E\left[\displaystyle\int_0^\infty\int_E \int_{\mathcal C_{[t-\sum_{i=1}^{n-1} S_j-s_n]}} h(\V_0,\Y^{(0)}_{[0,S_{1}[},\dots, s_n, \vv_n,\y^{(n)}_{[t-\sum_{i=1}^{n-1} S_j-s_n]})\1_{\left\{ {\sum_{j=1}^{n-1} S_j+s_n \leq t}\right\}}\right.\\
  & \left.\times \exp\left(-\displaystyle\int_0^{t-\sum_{j=1}^{n-1} S_j+s_n}\alpha(\y^{(n)}_u)d u\right)
L\left(\y^{(n)}_{[t-\sum_{i=1}^{n-1} S_j-s_n]}\right)dP^{\vv_n}_{[t-\sum_{i=1}^{n-1} S_j-s_n]}(\y^{(n)})\right.\\
 &  \left.\times\, k(\Y^{(n-1)}_{s_n},\vv_n)\nu(\Y^{(n-1)}_{s_n},d\vv_n)
\alpha(\Y^{(n-1)}_{s_n}) \exp\left(-\displaystyle\int_0^{s_n}\alpha(\Y_u^{(n-1)})d u\right)ds_n 
 \right].
\end{align*}
We can continue successively this process, using the (conditional) density of each move process $\Y^{(i)}$ on $[0,S_i]$, of each post-jump location $\V_i$, and of each inter-jump time $S_i$, and we obtain
\begin{align*}
 \E&  [h(\tilde\Xx)\1_{\left\{ {N_t=n}\right\}} ]\\
  =&  \int_{{\mathcal W}_n} h(\vv_0,\y^{(0)}_{[s_1]},\dots,s_n,\vv_n, \y^{(n)}_{[t-\sum_{i=1}^n s_j]})\1_{\left\{ {\sum_{j=1}^{n}s_j\leq t}\right\}}\\
  &\times \exp\left(-\displaystyle\int_0^{t-\sum_{j=1}^ns_j}\alpha(\y^{(n)}_u)d u\right)
      \times L(\y^{(n)}_{[t-\sum_{i=1}^n s_j]}) dP^{\vv_n}_{[t-\sum_{i=1}^n s_j]}(\y^{(n)})\\
      &\times k(\y^{(n-1)}_{s_n},\vv_n) \nu(\y^{(n-1)}_{s_n},d\vv_n) \alpha(\y^{(n-1)}_{s_n}) \exp\left(-\displaystyle\int_0^{s_n}\alpha(\y_u^{(n-1)})d u\right)ds_n \\
    & \times \prod_{i=1}^{n-1} \left[L(\y^{(i)}_{[s_{i+1}]})dP^{\vv_i}_{[s_{i+1}]}(\y^{(i)}) k(\y^{(i-1)}_{s_i},\vv_i) \nu(\y^{(i-1)}_{s_i},d\vv_i) \alpha(\y^{(i-1)}_{s_i})\right.\\
    &\left.\times \exp\left(-\displaystyle\int_0^{s_i}\alpha(\y_u^{(i-1)})d u\right)d s_i\right]
    L(\y^{(0)}_{[s_{1}]}) dP^{\vv_0}_{[s_1]}(\y^{(0)}) p(d\vv_0),
\end{align*}
where $p$ denotes the distribution of the initial state on $E$. 
 If we  come back to the initial expression \eqref{Xx} of $\Xx$ from that of $\tilde \Xx$, 
 using the shortcut $L(\x_{[t_{i},t_{i+1}[})$ for $L(\y^{(i)}_{[s_{i+1}]})=L(\z^{(i)}_{[0,t_{i+1}-t_{i}]},\m^{(i)})$ in view of $\x_t=\{(\z^{(i)}_{t-T_i},\m^{(i)})\}$ for all $t\in [t_{i},t_{i+1}[$, we obtain 
 \begin{align*}
& \E  [h(\Xx)\1_{\left\{ {N_t=n}\right\}} ]
=  \int_{{\mathcal W}_n} h(\x_0,\x_{[0,t_1[},\dots,t_n-t_{n-1},\x_{t_n},\x_{[t_n,t]})\1_{t_1\leq \dots\leq t_n\leq t} \\
&\times \exp\left(-\int_0^t  \alpha(\x_u)d u\right) L(\x_{[t_n,t]}) \prod_{i=1}^n  \left(L(\x_{[t_{i-1},t_i[}) k(\x_{t_i^-},\x_{t_i})\alpha(\x_{t_i^-})\right)\eta^{(n)}(d\xx),\end{align*}
where  $\eta^{(n)}$ is the measure given for 
 $\xx=(\x_0,\x_{[0,t_1[},\dots,t_n-t_{n-1},\x_{t_n},\x_{[t_n,t]})$ by
 \begin{multline*}\eta^{(n)}(d\xx)=\1_{t_1\leq \dots\leq t_n\leq t} dP^{\x_{t_n}}_{[t_n,t]}(\x)\nu(\x_{t_n^-},d\x_{t_n})dt_n \\
 \times \prod_{i=1}^{n-1} \left(dP^{\x_{t_i}}_{[t_i,t_{i+1}[}(\x)\nu(\x_{t_i^-},\x_{t_i})dt_i\right) dP^{\x_0}_{[0,t_1[}(\x)p(d\x_0).\end{multline*}
Here, just as we have used the natural notation $L(\x_{[t_{i},t_{i+1}[})$ for $L(\y^{(i)}_{[s_{i+1}]})$, we use $dP^{\x_{T_i}}_{[t_i,t_{i+1}[}(\x)$ for $dP^{\x_{T_i}}_{[s_{i+1}]}(\y^{(i)})$.

 Since for any $h$, $\E(h(\Xx))=\sum_{n\geq 0}\E(h(\Xx)\1_{\left\{ {N_t=n}\right\}}),$ we have proven that the likelihood of $\Xx$ given by \eqref{Xx}, with respect to the underlying measure
$ \eta(.)=\displaystyle\sum_{n\in \N}\eta^{(n)}(.\cap {\mathcal W}_n)$
 is given by 
$$\mathcal L(\Xx)=\exp\left(-\int_0^t  \alpha(\X_u)d u\right) L(\X_{[t_{N_t},t]}) \prod_{i=1}^{N_t}  \left(L(\X_{[T_{i-1},T_i[}) k(\X_{T_i^-},\X_{T_i})\alpha(\X_{T_i^-})\right).$$
If we replace  $\alpha$ and $k$ by their specific expression in \eqref{K_JM}, we obtain the result in the case of a BDMM process. Given that there is a one to one correspondance between the representation \eqref{Xx} and the trajectory $(\X_s)_{0\leq s\leq t}$, we can view the above likelihood as the likelihood of $(\X_s)_{0\leq s\leq t}$, where $\eta$ translates to a measure on the space of cadlag function from $[0,t]$ to $E$.

\section{Proof of Theorem \ref{LAN}}\label{proof LAN}

Let $t>0$ and $h\in\R^\ell$. We let $\vartheta_t= \vartheta+  h/ \sqrt{t}$ and we denote 
$$	\Lambda^{ \vartheta_{t}/ \vartheta}(t) = \log \frac {\mathcal{L}_{t}(\vartheta+\frac{h}{\sqrt t})}  {\mathcal{L}_{t}(\vartheta)}.$$
By Theorem \ref{densite} the log-likelihood ratio process is given by: 
\begin{multline}\label{big lambda}
 \Lambda^{ \vartheta_{t}/ \vartheta}(t)=\Lambda_{L}^{ \vartheta_{t}/ \vartheta}(t)+\Lambda_{ \beta}^{ \vartheta_{t}/ \vartheta}(t)+\Lambda_{\delta}^{ \vartheta_{t}/ \vartheta}(t)+\Lambda_{\tau}^{ \vartheta_{t}/ \vartheta}(t)\\+\Lambda_{k_{ \beta}}^{ \vartheta_{t}/ \vartheta}(t)+\Lambda_{k_{ \delta}}^{ \vartheta_{t}/ \vartheta}(t)+\Lambda_{k_{ \tau}}^{ \vartheta_{t}/ \vartheta}(t),
 \end{multline}
where, whether $\gamma$ is $\beta$, $\delta$ or $\tau$,
\begin{align} 
\Lambda_{ \gamma}^{ \vartheta_{t}/ \vartheta}(t)&= \int_0^t \log\left( \frac{ \gamma( \vartheta_{t},\X_{s^{-}})}{ \gamma( \vartheta,\X_{s^{-}})}\right)dN_s^\gamma - \int_{0}^t ( \gamma(\vartheta_{t},\X_{s}) - \gamma(\vartheta,\X_{s}) )ds, \label{lambda2}\\
\Lambda_{k_{ \gamma}}^{ \vartheta_{t}/ \vartheta}(t)&= \int_0^t  \log\left( \frac{  k_{ \gamma}\left( \vartheta_{t},\X_{s^{-}}, \X_{s}\right),
}{ k_{ \gamma}\left(\vartheta,\X_{s^{-}}, \X_{s}\right)}\right)dN_s^\gamma, \label{lambdakg}\\
\Lambda_{ L}^{ \vartheta_{t}/ \vartheta}(t)&=\log \frac{L(\vartheta_t,\X _{[T_{N_{t}},t]})}{L(\vartheta,\X _{[T_{N_{t}},t]})} +\sum_{i=0}^{N_{t}-1} \log \frac{L(\vartheta_t,\X_{[T_{i},T_{i+1}[})}{L(\vartheta,\X_{[T_{i},T_{i+1}[})}.\label{lambdaGamma0}
\end{align}

The proof of Theorem~\ref{LAN} follows a standard scheme as in \cite{luschgy} and \cite{evalan}. The main step consists in proving the LAN decomposition  for each of the three terms \eqref{lambda2}, \eqref{lambdakg} and \eqref{lambdaGamma0}  above, which is the statement of the three following lemmas. Then  \eqref{eqlambda} is deduced by combining these decompositions, as carried out next.

\begin{lem}\label{lem_gamma}
Under the same assumptions as in Theorem~\ref{LAN}, we have for any $t>0$, $h\in\R^\ell$ and $\vartheta\in\Theta$,
\begin{align}\label{decompogamma}
\Lambda_{ \gamma}^{ \vartheta_{t}/ \vartheta}(t)=h^T M_{ \gamma}^{ \vartheta}(t)- \frac{ 1}{ 2}h^T\langle M_{ \gamma}^{ \vartheta} \rangle(t) h+Rem^ \gamma_{ \vartheta,h}(t),
\end{align}
where $Rem^\gamma_{\vartheta,h}(t)$ tends to $0$ in $\P_{\x_0}^\vartheta$-probability when $t \to +\infty$, and   $(M_{ \gamma}^{\vartheta}(t))_{t>0}$ is a martingale with respect to $\mathcal F_{t}$ given by \eqref{M_gamma}. 
\end{lem}

\begin{lem}\label{lem_kgamma}
Under the same assumptions as in Theorem~\ref{LAN}, we have for any $t>0$, $h\in\R^\ell$ and $\vartheta\in\Theta$,
\begin{align}\label{decompokgamma}
\Lambda_{k_{ \gamma} }^{ \vartheta_{t}/ \vartheta}(t)=h^T M_{ k_{ \gamma}}^{ \vartheta}(t)- \frac{ 1}{ 2}h^T\langle M_{ k_{\gamma}}^{ \vartheta} \rangle(t) h+Rem^ {k_{\gamma}}_{ \vartheta,h}(t),
\end{align}
where $Rem^{k_\gamma}_{\vartheta,h}(t)$ tends to $0$ in $\P_{\x_0}^\vartheta$-probability when $t \to +\infty$, and $(M_{ k_\gamma}^{ \vartheta}(t))_{t>0}$ is a martingale with respect to $\mathcal F_{t}$ given by \eqref{M_kgamma}. 
\end{lem}

\begin{lem}\label{lem_L}
Under the same assumptions as in Theorem~\ref{LAN}, we have for any $t>0$, $h\in\R^\ell$ and $\vartheta\in\Theta$,
\begin{align}\label{decompokgamma}
\Lambda_{L }^{ \vartheta_{t}/ \vartheta}(t)=h^T M_{ L}^{\vartheta}(t)- \frac{ 1}{ 2}h^T\langle M_{ L}^{\vartheta} \rangle(t) h+Rem^{L}_{ \vartheta,h}(t),
\end{align}
where $Rem^{L}_{\vartheta,h}(t)$ tends to $0$ in $\P_{\x_0}^\vartheta$-probability when $t \to +\infty$, and  $(M_{ L}^{ \vartheta}(t))_{t>0}$ is a martingale with respect to $\mathcal F_{t}$ given by \eqref{M_L}.
\end{lem}

The proofs of these three lemmas are postponed to the end of this  section.  Let us deduce \eqref{eqlambda}.  From their statements and \eqref{big lambda}, we have 
$$\Lambda^{ \vartheta_{t}/ \vartheta}(t) = h^T M^{ \vartheta}(t) - \frac{ 1}{ 2}h^T\sum_{i\in\mathcal S} \langle M_{ i}^{\vartheta} \rangle(t) h+Rem_{\vartheta,h}(t),$$
where $\mathcal S=\{L,\beta,\delta,\tau,k_\beta,k_\delta,k_\tau\}$, $M^{ \vartheta}(t) =\sum_{i\in\mathcal S} M_i^{ \vartheta}$ and $Rem_{ \vartheta,h}(t)=\sum_{i\in\mathcal S}  Rem^i_{\vartheta,h}(t)$. By Lemmas~\ref{lem_gamma}, \ref{lem_kgamma} and \ref{lem_L}, we know that 
$Rem_{\vartheta,h}(t)$ tends to $0$ in $\P_{\x_0}^\vartheta$-probability. To show the decomposition \eqref{eqlambda}, it remains to prove that 
$$\langle M^{\vartheta}\rangle(t) = \sum_{i\in\mathcal S} \langle M_{ i}^{\vartheta} \rangle(t),$$
which, given the definition of $M^{\vartheta}(t)$, boils down to proving that
\begin{equation}\label{cross angle 0}\forall i,j \in \mathcal{S},\: i\neq j,\quad  \langle M_{ i}^{ \vartheta} , M_{ j}^{ \vartheta} \rangle(t) =0.\end{equation}

First, for $i=L$, since $M_{L}^{\vartheta}(t)$ is a continuous martingale and for $j\neq L$, $M_{ j}^{ \vartheta}$ has finite variations, we have that $\langle M_{L}^{ \vartheta} , M_{ j}^{ \vartheta} \rangle(t)=0$ for any $j\neq L$, see \cite[Chapter 8.9]{klebaner}.  Second, for all $i,j\in \mathcal S\setminus\{L\}$, $i\neq j$, we have
$$ [M_{i}^{ \vartheta}, M_{j}^{ \vartheta}](t)=\displaystyle\sum_{s\leq t} \Delta M_{i}^{ \vartheta}(s) \Delta M_{j}^{ \vartheta}(s),$$
where $\Delta M_{i}^{ \vartheta}(s)=M_{i}^{ \vartheta}(s)-M_{i}^{ \vartheta}(s^-)$. 
For $\gamma\neq\gamma'$, $M_{ \gamma}^{ \vartheta}(t)$ and $ M_{ \gamma'}^{ \vartheta}(t)$ do not have any jump in common, and similarly for  $M_{ \gamma}^{ \vartheta}(t)$ and $M_{ k_{\gamma'}}^{ \vartheta}(t)$, since a birth, a death or a mutation cannot occur at the same moment almost surely. Consequently $[M_{\gamma}^{ \vartheta}, M_{\gamma'}^{ \vartheta}](t)=[M_{\gamma}^{ \vartheta}, M_{k_\gamma'}^{ \vartheta}](t)=0$ for any $\gamma\neq\gamma'$. This implies \eqref{cross angle 0} for $i=\gamma$ and $j\in\{\gamma',k_\gamma'\}$, $\gamma\neq\gamma'$.  To complete the proof of \eqref{cross angle 0}, it remains to show that $\langle M_{ \gamma}^{ \vartheta}, M_{k_{ \gamma}}^{ \vartheta}\rangle (t) =0$. We have that 
\begin{align*}
[M_{ \gamma}^{ \vartheta}, M_{k_{ \gamma}}^{ \vartheta}] (t)=\frac{1}{ t} \int_{0}^{t}  \frac{ \nabla \gamma( \vartheta,\X_{s^-})}{ \gamma( \vartheta,\X_{s^-})} \frac{\left( \nabla k_\gamma ( \vartheta,\X_{s^-},\X_s)\right)^T}{ k_\gamma ( \vartheta,\X_{s^-},\X_s)} dN^{\gamma}_{s},
\end{align*}
so that by Lemma~\ref{martingale_double}, the associated angle bracket (corresponding  to the compensator of the square bracket, see  \cite{klebaner})  is given by 
\begin{align*}
\langle M_{ \gamma}^{ \vartheta}, M_{k_{ \gamma}}^{ \vartheta}\rangle (t)&=\frac{1}{ t} \displaystyle \int_{0}^{t} \nabla \gamma( \vartheta,\X_{s}) \int_{E} \left(\nabla k_\gamma ( \vartheta,\X_{s},\y)\right)^T\nu_{ \gamma}(\X_{s},d\y)ds.
\end{align*}
Remember that  for any $\vartheta\in\Theta$ and any $\x\in E$, $\int_E k_\gamma(\vartheta,\x,\y)\nu(\x,d\y)=1$. This implies $ \int_{E} \nabla k_\gamma ( \vartheta,\x,\y)\nu_{ \gamma}(\x,d\y)=0$ if interchange of integration and differentation is valid. In particular, it can be showed that the latter holds under Condition~\ref{cond:Btilde}. Consequently, $\langle M_{ \gamma}^{ \vartheta}, M_{k_{ \gamma}}^{ \vartheta}\rangle (t) =0$ and \eqref{cross angle 0} is proven, which completes the proof of \eqref{eqlambda}.

The last statement of Theorem~\ref{LAN}, that is the convergence \eqref{thmeq2}, is a direct consequence of Theorems 4.12 and 4.22 of \cite{hopfnereva}, thanks to our hypothesis \ref{hypH} that implies that $\textnormal{\O}$ is a recurrent atom for $(\X_t)_{t>0}$.

\subsection{Proof of Lemma~\ref{lem_gamma}}

Let us prove that $Rem^ \gamma_{ \vartheta,h}(t)$ converges to 0 in $\P_{\x_0}^\vartheta$-probability, where
$$Rem^ \gamma_{ \vartheta,h}(t) = \Lambda_{ \gamma}^{ \vartheta_{t}/ \vartheta}(t) - h^T M_{ \gamma}^{ \vartheta}(t) + \frac{ 1}{ 2}h^T\langle M_{ \gamma}^{ \vartheta} \rangle(t) h,$$
and $M_{ \gamma}^{ \vartheta}(t)$ is given by  \eqref{M_gamma}. From \eqref{lambda2}, we get 
\begin{multline*}
\Lambda_{ \gamma}^{ \vartheta_{t}/ \vartheta}(t)= \int_{0}^{t} \left[  \log\left(\frac{ \gamma( \vartheta_{t},\X_{s^-})}{ \gamma( \vartheta,\X_{s^-})}\right)-\frac{ \gamma( \vartheta_{t},\X_{s^-})}{ \gamma( \vartheta,\X_{s^-})}+1\right] dN^{\gamma}_{s} \\ + \int_{0}^{t} \left[ \frac{ \gamma( \vartheta_{t},\X_{s^-})}{ \gamma( \vartheta,\X_{s^-})}-1\right]d\tilde N^{\gamma}_{s},
\end{multline*}
where  $\tilde N^{\gamma}_{t}=N^{\gamma}_{t}- \displaystyle\int_{0}^t \gamma(\vartheta, \X_{s})ds$ is a $\mathcal F_t$-martingale. 
We have in particular $[\tilde N^{\gamma}]_t=N_t^\gamma$ and $\langle \tilde N^{\gamma}\rangle_t=\int_0^t \gamma(\vartheta, \X_{s})ds$, see for instance \cite{klebaner}.   We deduce that
\begin{align*}
[M_{ \gamma}^{ \vartheta}](t)&=\frac{ 1}{ t}\int_{0}^{t}\frac{ (\nabla \gamma) (\nabla \gamma)^T}{ \gamma^2}( \vartheta,\X_{s^-})d[\tilde N^{\gamma}]_{s} \\ & =\frac{ 1}{ t}\int_{0}^{t}\frac{ (\nabla \gamma) (\nabla \gamma)^T}{ \gamma^2}( \vartheta,\X_{s^-})dN^\gamma_{s}.
\end{align*}
Note that $\gamma(\vartheta,\X_{T_i^-})>0$ for all $\gamma$-jump time $T_i$, so that the previous formula remains true with the addition of $\1_{\gamma(\vartheta,\X_{s^-})>0}$ in the integrand. So we get 
\begin{align*}
\langle M_{ \gamma}^{ \vartheta} \rangle(t)&
=\frac{ 1}{ t}\int_{0}^{t}\1_{\gamma>0}\frac{ (\nabla \gamma) (\nabla \gamma)^T}{ \gamma}( \vartheta,\X_{s})ds,
\end{align*}
as claimed in \eqref{brackets_gamma}, which well satisfies the property that $[M_{ \gamma}^{ \vartheta}](t)-\langle M_{ \gamma}^{ \vartheta} \rangle(t)$ is a martingale. 
Using these representations, we may write 
\begin{align*}
&Rem^ \gamma_{ \vartheta,h}(t) \\ &= \Lambda_{ \gamma}^{ \vartheta_{t}/ \vartheta}(t) - h^T M_{ \gamma}^{ \vartheta}(t) + \frac{ 1}{ 2}h^T\left(\langle M_{ \gamma}^{ \vartheta} \rangle(t) -  [ M_{ \gamma}^{ \vartheta} ](t)\right)h  + \frac{ 1}{ 2}h^T[M_{ \gamma}^{ \vartheta}](t) h
\\
&= S_t + A_t + R_t,
\end{align*}
where
\begin{align*}
S_{t}&=\int_{0}^{t} \left( \frac{ \gamma( \vartheta_{t},\X_{s^-})}{ \gamma( \vartheta,\X_{s^-})}-1-\frac{ h^T}{ \sqrt{t}}\frac{ \nabla \gamma( \vartheta,\X_{s^-})}{ \gamma( \vartheta,\X_{s^-})}\right)d\tilde N^{\gamma}_{s},\\
A_{t}&=\frac{ 1}{ 2t}\int_{0}^{t}\1_{\gamma>0}\frac{ (h^T\nabla \gamma)^2 }{ \gamma}( \vartheta,\X_{s})ds -\frac{ 1}{2 t}\int_{0}^{t} \frac{ (h^T\nabla \gamma)^2 }{ \gamma^2}( \vartheta,\X_{s^-})dN^{\gamma}_{s},\\
R_{t}&= \int_{0}^{t}\left(  \log\left(\frac{ \gamma( \vartheta_{t},\X_{s^-})}{ \gamma( \vartheta,\X_{s^-})}\right)-\frac{ \gamma( \vartheta_{t},\X_{s^-})}{ \gamma( \vartheta,\X_{s^-})}+1+\frac 1{2t} \frac{ (h^T\nabla \gamma)^2 }{ \gamma^2}( \vartheta,\X_{s^-})\right)dN^{\gamma}_{s}.
\end{align*}
To complete the proof, we show below that each of these three terms tends to 0 in $\P_{\x_0}^\vartheta$-probability.

\bigskip

\noindent\underline{First step:} \textit{Convergence of $S_{t}$ to $0$ in probability.} \\
By the Markov inequality, we have for any $a>0$, 
$$\mathbb{P}_{\x_0}^\vartheta\left( |S_{t}|>a\right)\leq \frac{ 1}{ a^2} \mathbb{E}_{\x_0}^\vartheta\left((S_{t})^2\right)= \frac{ 1}{ a^2} \mathbb{E}_{\x_0}^\vartheta\left( \langle S\rangle_{t}\right),$$
where the last equality comes from the fact that $(S_{t})^2-\langle S\rangle_{t}$ is by definition a centred martingale (see \cite{klebaner}).  We prove in the following that $\mathbb{E}_{\x_0}^\vartheta\left( \langle S\rangle_{t}\right)$ tends to 0. We have
\begin{align*}
[S]_{t}&= \int_{0}^{t} \left( \frac{ \gamma( \vartheta_{t},\X_{s^-})}{ \gamma( \vartheta,\X_{s^-})}-1-\frac{ h^T}{ \sqrt{t}}\frac{ \nabla \gamma( \vartheta,\X_{s^-})}{ \gamma( \vartheta,\X_{s^-})}\right)^2 dN^\gamma_{s}\\
&=\int_{0}^{t} \left( \frac{ \gamma( \vartheta_{t},\X_{s^-})}{ \gamma( \vartheta,\X_{s^-})}-1-\frac{ h^T}{ \sqrt{t}}\frac{ \nabla \gamma( \vartheta,\X_{s^-})}{ \gamma( \vartheta,\X_{s^-})}\right)^2 \1_{ \gamma( \vartheta,\X_{s^-})>0}dN^\gamma_{s},
\end{align*}
so that
\begin{align}\label{Yn}
\langle S\rangle _{t}
&= \int_{0}^{t} \left( \frac{ \gamma( \vartheta_{t},\X_{s})}{ \gamma( \vartheta,\X_{s})}-1-\frac{ h^T}{ \sqrt{t}}\frac{ \nabla \gamma( \vartheta,\X_{s})}{ \gamma( \vartheta,\X_{s})}\right)^2  \1_{ \gamma( \vartheta,\X_{s})>0} \gamma( \vartheta, \X_{s})ds.
\end{align}
Note that for any $\vartheta,\vartheta'\in\Theta$ and any $\x$,
\begin{align*}
\gamma( \vartheta',\x)- \gamma( \vartheta,\x)=\int_{0}^1 \left(\nabla  \gamma( \vartheta+t( \vartheta'- \vartheta),\x)\right)^T( \vartheta'- \vartheta) dt.
\end{align*}
Thus we have for all $\x\in E$ and $\vartheta, \vartheta' \in \Theta$ such that $\gamma( \vartheta,\x)>0$:
\begin{align}\label{dl gamma}
\bigg( \frac{ \gamma( \vartheta',\x)}{ \gamma( \vartheta,\x)}-1- & \frac{\left(\nabla  \gamma( \vartheta,\x)\right)^T( \vartheta'- \vartheta)}{ \gamma(\vartheta,\x)}\bigg)^2 \nonumber \\
&=  \left( \displaystyle \int_{0}^1 \frac{ \left(\nabla  \gamma( \vartheta+t( \vartheta'- \vartheta),\x)-\nabla  \gamma( \vartheta,\x)\right)^T( \vartheta'- \vartheta)}{ \gamma( \vartheta,\x)}dt\right)^2\nonumber\\
&\leq  \left\| \vartheta'- \vartheta \right\|^2 f^ \gamma_{ \vartheta}(\x, \left\| \vartheta'- \vartheta \right\|),
\end{align}
where $f^ \gamma_{ \vartheta}(\x, \rho)$ is given by \eqref{f_gamma}.
Applying this inequality in \eqref{Yn} with $\vartheta'=\vartheta_t= \vartheta+  h/ \sqrt{t}$, we obtain
\begin{align}\label{ineqf}
\langle S\rangle _{t}&\leq \frac{ \|h\| ^2}{ t} \int_{0}^{t} f^ \gamma_{ \vartheta}\left(\X_{s},\frac{ \|h\|}{ \sqrt{t}}\right) \gamma( \vartheta,\X_{s}) ds.
\end{align}
Since  $f^ \gamma_{ \vartheta}(\x, .)$ is increasing,  we have by condition~\ref{cond:B} that for any $\vartheta\in\Theta$ and $t$ large enough, 
\begin{equation}\label{borneinfini}
\sup_{\x\in E} \left |f^ \gamma_{ \vartheta}\left(\x,\frac{ \|h\|}{ \sqrt{t}}\right) \gamma( \vartheta,\x)\right | < d(\vartheta),
\end{equation}
for some $d(\vartheta)>0$. We can thus apply Corollary~\ref{corH} in Appendix~\ref{appendix:ergodicity} to the right-hand side of \eqref{ineqf}, thanks to the ergodic assumption \ref{hypH},   to get
\begin{align*}
\mathbb{E}_{\x_{0}}^\vartheta\left( \left|\langle S\rangle _{t} \right|\right) 
&\leq \frac{ \|h\|^2 c(\x_0)d(\vartheta)}{ t}+  \|h\|^2  \int_{E}\left| f^ \gamma_{ \vartheta}\left(\x,\frac{ \|h\|}{ \sqrt{t}}\right) \gamma( \vartheta,\x)\right| \mu_{\infty}^\vartheta(d\x).\end{align*}
In view of \eqref{borneinfini}, and since for any $\x\in E$,  $f^ \gamma_{ \vartheta}(\x, \|h\|/ \sqrt{t})\to 0$ when $t\to 0$, we can conclude by the dominated convergence theorem.

\bigskip

\noindent\underline{Second step:} \textit{Convergence of $A_{t}$ to $0$ in probability.} \\
Let $ \epsilon>0$. We have $2A_{t} = A_1 -  A_2 - A_3$ where
\begin{align}
A_1 & = \frac{ 1}{t}\int_{0}^{t}\1_{\gamma>0} \frac{ (h^T\nabla \gamma)^2 }{ \gamma}( \vartheta,\X_{s})\1_{\left\{ \left| \frac{h^T\nabla \gamma}{\gamma}( \vartheta,\X_{s})\right| > \epsilon\sqrt{t}\right\}}ds,\nonumber\\
A_2 & =  \frac{ 1}{ t}\int_{0}^{t} \frac{ (h^T\nabla \gamma)^2 }{ \gamma^2}( \vartheta,\X_{s^-})\1_{\left\{ \left| \frac{h^T\nabla \gamma}{\gamma}( \vartheta,\X_{s^-})\right| > \epsilon\sqrt{t}\right\}}dN^{\gamma}_{s},\label{A2}\\
A_3 & =   \frac{ 1}{ t} \int_{0}^{t}\1_{\gamma>0} \frac{ (h^T\nabla \gamma)^2 }{ \gamma^2}( \vartheta,\X_{s^-})\1_{\left\{ \left| \frac{h^T\nabla \gamma}{\gamma}( \vartheta,\X_{s^-})\right| \leq \epsilon\sqrt{t}\right\}} d\tilde N^{\gamma}_s.\nonumber
\end{align}
Note that $\E_{\x_0}^\vartheta(|A_1|)=\E_{\x_0}^\vartheta(A_1)=\E_{\x_0}^\vartheta(A_2)=\E_{\x_0}^\vartheta(|A_2|)$. We prove below that for any $\epsilon>0$, $\E_{\x_0}^\vartheta(|A_1|)\to 0$ as $t\to \infty$ and then we deal with the convergence in probability of $A_3$. 

By  condition \ref{cond:A}, we know that given $\vartheta\in\Theta$,  the integrand in $A_1$ is uniformly bounded in $\X_s$ by some $d(\vartheta)>0$. Thanks to \ref{hypH}  and Corollary~\ref{corH}, we get
$$\E_{\x_0}^\vartheta(|A_1|)\leq \frac{c(\x_0)d(\vartheta)}{t} +  \int_E \left| \1_{\gamma>0} \frac{ (h^T\nabla \gamma)^2 }{ \gamma}( \vartheta,\x)\1_{\left\{ \left| \frac{h^T\nabla \gamma}{\gamma}( \vartheta,\x)\right| > \epsilon\sqrt{t}\right\}}\right| \mu_{\infty}^\vartheta(d\x).$$
The latter integrand is bounded thanks to condition \ref{cond:A} and converges to 0 as $t\to\infty$ for any given $\x$ and $\vartheta$. We thus conclude by the dominated convergence theorem that $\E_{\x_0}^\vartheta(|A_1|)\to 0$ as $t\to\infty$. 

Concerning $A_3$, since it is a centred martingale, we have for any $a>0$, 
$$\P_{\x_0}^\vartheta(|A_3|>a)\leq \frac 1 {a^2} \E_{\x_0}^\vartheta (|A_3|^2) =\frac 1 {a^2} \E_{\x_0}^\vartheta (\langle A_3\rangle)$$ 
where 
\begin{align*}
\langle A_3\rangle &=  \frac{ 1}{ t^2} \int_{0}^{t}\1_{\gamma>0} \frac{ (h^T\nabla \gamma)^4 }{ \gamma^3}( \vartheta,\X_{s})\1_{\left\{ \left| \frac{h^T\nabla \gamma}{\gamma}( \vartheta,\X_{s})\right| \leq \epsilon\sqrt{t}\right\}} ds\\
&\leq \frac{\epsilon^2}{t} \int_{0}^{t}\1_{\gamma>0}  \frac{ (h^T\nabla \gamma)^2 }{ \gamma}( \vartheta,\X_{s}) ds.
\end{align*}
We can again apply  Corollary~\ref{corH} to the latest integral, thanks to conditions \ref{hypH} and \ref{cond:A}, and we obtain 
$$ \E_{\x_0}^\vartheta (\langle A_3\rangle) \leq \epsilon^2 C(\x_0,\vartheta,h)$$
where $C(\x_0,\vartheta,h)>0$. So for any $a>0$ and any $\eta>0$, we can choose $\epsilon$ so that $\P_{\x_0}^\vartheta(|A_3|>a)\leq \eta$. For these choices, we can further choose $t$ large enough so that  $\P_{\x_0}^\vartheta(|A_1|>a)\leq \eta$ and $\P_{\x_0}^\vartheta(|A_2|>a)\leq \eta$, because for any $\epsilon>0$, $\E_{\x_0}^\vartheta(|A_1|)=\E_{\x_0}^\vartheta(|A_2|)\to 0$ as $t\to \infty$. This entails the convergence of $A_{t}$ to $0$ in $\P_{\x_0}^\vartheta$-probability.

\bigskip

\noindent\underline{Third step:} \textit{Convergence of $R_{t}$ to $0$ in probability.} \\
Denoting by $\varphi:x>-1\mapsto \log(1+x)-x+ x^2/2$, we have $R_{t}=R_{1}-R_{2}/2$ 
with
\begin{align}
R_{1}&=\int_{0}^{t}  \varphi \left( \frac{ \gamma( \vartheta_{t},\X_{s^-})}{ \gamma(\vartheta,\X_{s^-})}-1\right)dN^{\gamma}_{s},\nonumber\\
R_{2}&= \int_{0}^{t} \left(\left[\frac{ \gamma( \vartheta_{t},\X_{s^-})}{ \gamma(\vartheta,\X_{s^-})}-1\right]^2 -\left[ \frac{ h^T}{ \sqrt{t}}\frac{ \nabla \gamma( \vartheta,\X_{s^-})}{ \gamma( \vartheta,\X_{s^-})}\right]^2\right) dN^{\gamma}_{s}.\label{R2}
\end{align}
We start by proving that $R_2$ tends to 0 in probability. Write $R_2=\int (f^2-g^2)dN^\gamma_s$ with obvious notations. We have by the Cauchy-Schwartz inequality
and the fact that  $(f+g)^2=(f-g+2g)^2\leq 2(f-g)^2+8g^2$, 
\begin{align*}
R_2^2 & \leq \int (f-g)^2dN^\gamma_s \int (f+g)^2dN^\gamma_s \\
&\leq  \left(\int (f-g)^2dN^\gamma_s\right)\left(8 \int (f-g)^2dN^\gamma_s + 8 \int g^2 dN^\gamma_s\right).
\end{align*}
To prove  the convergence of $R_2$, it is then sufficient to show that $\int (f-g)^2dN^\gamma_s = o_{\P_{\x_0}^\vartheta}(1)$ and $ \int g^2 dN^\gamma_s=O_{\P_{\x_0}^\vartheta}(1)$. We have
\begin{align}\label{f-g}
\E_{\x_0}^\vartheta & \int (f-g)^2dN^\gamma_s \nonumber \\
& =  \E_{\x_0}^\vartheta \int_{0}^{t} \left(\frac{ \gamma( \vartheta_{t},\X_{s^-})}{ \gamma(\vartheta,\X_{s^-})}-1 - \frac{ h^T}{ \sqrt{t}}\frac{ \nabla \gamma( \vartheta,\X_{s^-})}{ \gamma( \vartheta,\X_{s^-})}\right)^2 dN^{\gamma}_{s} \\
& =\E_{\x_0}^\vartheta   \int_{0}^{t} \left( \frac{ \gamma( \vartheta_{t},\X_{s})}{ \gamma( \vartheta,\X_{s})}-1-\frac{ h^T}{ \sqrt{t}}\frac{ \nabla \gamma( \vartheta,\X_{s})}{ \gamma( \vartheta,\X_{s})}\right)^2  \1_{ \gamma( \vartheta,\X_{s})>0} \gamma( \vartheta, \X_{s})ds, \nonumber
 \end{align}
which is exactly $\E_{\x_0}^\vartheta(\langle S\rangle _{t})$, see \eqref{Yn}. We have already proven that this term tends to 0, so $\int (f-g)^2dN^\gamma_s = o_{\P_{\x_0}^\vartheta}(1)$. Now, for any $\eta>0$, by the Markov's inequality, we can choose $M>0$ such that
\begin{align*}
\P_{\x_0}^\vartheta \left(\int g^2 dN^\gamma_s>M\right) \leq \frac 1M \E_{\x_0}^\vartheta \left(\frac 1t\int_0^{t} \frac{  (h^T\nabla \gamma( \vartheta,\X_{s}))^2}{ \gamma( \vartheta,\X_{s})}ds\right)\leq \frac 1 M C(\vartheta,h),
\end{align*}
where $C(\vartheta,h)$ is a positive upper-bound deduced from condition~\ref{cond:A}. This proves that $ \int g^2 dN^\gamma_s=O_{\P_{\x_0}^\vartheta}(1)$ and completes the proof that  $R_2$ tends to 0 in probability. 

It remains to address the convergence of $R_1$. For any $\epsilon>0$, we have $R_1=R_{11}+R_{12}$ with
\begin{align}
R_{11}&=\int_{0}^{t} \varphi \left(\frac{ \gamma( \vartheta_{t},\X_{s^-})}{ \gamma(\vartheta,\X_{s^-})}-1\right)\1_{\left\{ \left|\frac{ \gamma( \vartheta_{t},\X_{s^-})}{ \gamma(\vartheta,\X_{s^-})}-1\right|> \epsilon\right\}}dN^{\gamma}_{s}, \\
R_{12}&= \int_{0}^{t} \varphi \left(\frac{ \gamma( \vartheta_{t},\X_{s^-})}{ \gamma(\vartheta,\X_{s^-})}-1\right)\1_{\left\{ \left|\frac{ \gamma( \vartheta_{t},\X_{s^-})}{ \gamma(\vartheta,\X_{s^-})}-1\right|\leq  \epsilon\right\}}dN^{\gamma}_{s}.\label{R12}
\end{align}
For $R_{11}$, note that for any $a>0$,
\begin{equation}\label{ab}
\P_{\x_0}^\vartheta (R_{11}>a) \leq  \P_{\x_0}^\vartheta  \left( \int_{0}^{t}\1_{\left\{ \left|\frac{ \gamma( \vartheta_{t},\X_{s^-})}{ \gamma(\vartheta,\X_{s^-})}-1\right|> \epsilon\right\}}dN^{\gamma}_{s} >b\right),
\end{equation}
for any $b<1$. This is because when $b<1$, 
\begin{multline*}
\left( \int_{0}^{t}\1_{\left\{ \left|\frac{ \gamma( \vartheta_{t},\X_{s^-})}{ \gamma(\vartheta,\X_{s^-})}-1\right|> \epsilon\right\}}dN^{\gamma}_{s} \leq b \right)\\ 
\Longrightarrow \left(\text{for all $\gamma$-jump time $T_i\leq t$ },\1_{\left\{ \left|\frac{ \gamma( \vartheta_{t},\X_{T_i^-})}{ \gamma(\vartheta,\X_{T_i^-})}-1\right|> \epsilon\right\}} =0 \right)
\end{multline*}
which implies that $R_{11}=0\leq a$. Moreover
\begin{align*}
\1_{\left\{ \left|\frac{ \gamma( \vartheta_{t},\x)}{ \gamma(\vartheta,\x)}-1\right|> \epsilon\right\}}
 &\leq\1_{\left\{ \left|\frac{ \gamma( \vartheta_{t},\x)}{ \gamma(\vartheta,\x)}-1-\frac{ h^T}{ \sqrt{t}}\frac{\nabla \gamma( \vartheta,\x)}{\gamma( \vartheta,\x)}\right|> \frac{ \epsilon}{ 2}\right\}} +\1_{\left\{ \left|\frac{ h^T}{ \sqrt{t}}\frac{\nabla \gamma( \vartheta,\x)}{\gamma( \vartheta,\x)}\right|> \frac{ \epsilon}{ 2}\right\}}.
 \end{align*}
Using this and the fact that for any $z$, $\1_{\left\{ |z|> \frac{ \epsilon}{ 2}\right\}}\leq 4 \epsilon^{-2}z^2\1_{\left\{ |z|> \frac{ \epsilon}{ 2}\right\}}$, we get
\begin{multline*}
\int_{0}^{t}\1_{\left\{ \left|\frac{ \gamma( \vartheta_{t},\X_{s^-})}{ \gamma(\vartheta,\X_{s^-})}-1\right|> \epsilon\right\}}dN^{\gamma}_{s}\\
\leq4 \epsilon^{-2}\int_{0}^{t}\left(\frac{ \gamma( \vartheta_{t},\X_{s^-})}{ \gamma(\vartheta,\X_{s^-})}-1-\frac{ h^T}{ \sqrt{t}}\frac{\nabla \gamma( \vartheta,\X_{s^-})}{\gamma( \vartheta,\X_{s^-})}\right)^2dN^{\gamma}_{s}\\
 +4 \epsilon^{-2}\int_{0}^{t} \left(\frac{ h^T}{ \sqrt{t}}\frac{\nabla \gamma( \vartheta,\X_{s^-})}{\gamma( \vartheta,\X_{s^-})}\right)^2\1_{\left\{ \left|\frac{ h^T}{ \sqrt{t}}\frac{\nabla \gamma( \vartheta,\X_{s^-})}{\gamma( \vartheta,\X_{s^-})}\right|> \frac{ \epsilon}{ 2}\right\}}dN^{\gamma}_{s}.
\end{multline*}
The first term is exactly $4 \epsilon^{-2}\int (f-g)^2dN^\gamma_s$, as already studied in \eqref{f-g}, and the second term is  $4 \epsilon^{-2}A_2$ (up to $\epsilon/2$ instead of $\epsilon$), see \eqref{A2}. Both terms tend to 0 in $\P_{\x_0}^\vartheta$-probability whatever the value of $\epsilon>0$. This proves that the right-hand side term in  \eqref{ab} tends to 0 for any $b<1$, which yields the convergence of $R_{11}$ in $\P_{\x_0}^\vartheta$-probability to 0. 

For $R_{12}$ given by \eqref{R12}, we choose $ \epsilon< \frac{ 1}{ 2}$ and since $|\varphi(x)|\leq 2|x|^3$ for $|x|\leq \frac{ 1}{ 2}$, we have
\begin{align*}
R_{12} &\leq \int_{0}^{t} 2\left|\frac{ \gamma( \vartheta_{t},\X_{s^-})}{ \gamma(\vartheta,\X_{s^-})} -1\right|^3\1_{\left\{ \left|\frac{ \gamma( \vartheta_{t},\X_{s^-})}{ \gamma(\vartheta,\X_{s^-})}-1\right|\leq \epsilon\right\}}dN^{\gamma}_{s}\\
&\leq2 \epsilon\displaystyle\int_{0}^{t} \left(\frac{ \gamma( \vartheta_{t},\X_{s^-})}{ \gamma(\vartheta,\X_{s^-})} -1\right)^2dN^{\gamma}_{s}.
\end{align*}
This last term is $2\epsilon(R_2 + \int g^2 dN^\gamma_s)$ if we use as previously the notation $R_2=\int (f^2-g^2)dN^\gamma_s$ for $R_2$ given by \eqref{R2}. We already know that $U:=2(R_2 + \int g^2 dN^\gamma_s)$ is a $O_{\P_{\x_0}^\vartheta}(1)$ as $t\to\infty$. This means that for any $\eta>0$, there exists $M$ such that for $t$ sufficiently large $\P_{\x_0}^\vartheta(U>M)<\eta$.  This implies that for any $a>0$ and $\eta>0$, we can choose $\epsilon$ small enough ($\epsilon<a/M$) so that $\P_{\x_0}^\vartheta(\epsilon U>a) <P_{\x_0}(U>M)<\eta$. So for any $a>0$ and $\eta>0$, we can choose $\epsilon$ so that $\P_{\x_0}^\vartheta(R_{12}>a) <\eta$ for $t$ large enough. The same inequality holds true for $R_{11}$ since we have proven that $R_{11}$ tends to 0 in probability for any $\epsilon>0$. So  $R_{1}=R_{11}+R_{12}$ tends to 0  in $\P_{\x_0}^\vartheta$-probability, which concludes the proof.

\subsection{Proof of Lemma~\ref{lem_kgamma}}
The proof follows the same scheme as the proof of Lemma~\ref{lem_gamma}. Consider 
$$Rem^ {k_{\gamma}}_{ \vartheta,h}(t)=\Lambda_{k_{ \gamma} }^{ \vartheta_{t}/ \vartheta}(t) - h^T M_{ k_{ \gamma}}^{ \vartheta}(t) + \frac{ 1}{ 2}h^T\langle M_{ k_{\gamma}}^{ \vartheta} \rangle(t) h,
$$
where $M_{ k_{ \gamma}}^{ \vartheta}(t)$ is given by \eqref{M_kgamma}. The fact that $M_{ k_{ \gamma}}^{ \vartheta}$ is a $\mathcal F_{t}$-martingale is a consequence of Lemma~\ref{martingale_double}. We deduce that
\begin{align*}
[M_{ k_{ \gamma}}^{ \vartheta}](t)&=\frac{ 1}{ t}\int_{0}^{t} \frac{ \nabla k_\gamma (\nabla k_\gamma)^T}{ k_\gamma^2}( \vartheta,\X_{s^-},\X_s) dN^{\gamma}_{s},\\
\langle M_{ k_{ \gamma}}^{ \vartheta}\rangle(t)&=\frac{1}{t}\int_{0} ^{t} \int_{E} \1_{k_\gamma>0} \frac{ \nabla k_\gamma (\nabla k_\gamma)^T}{ k_\gamma}( \vartheta,\X_{s},\y)\nu_{ \gamma}(\X_{s},d\y) \gamma( \vartheta, \X_{s})ds.
\end{align*}

Starting from  \eqref{lambdakg} and using the above representations we can decompose $Rem^ {k_{\gamma}}_{ \vartheta,h}(t)$ as in the proof of Lemma~\ref{lem_gamma}, that is $Rem^{k_\gamma}_{ \vartheta,h}(t)= S_t + A_t + R_t $, where
\begin{multline*}
S_{t}=\int_{0}^{t} \left( \frac{ k_{ \gamma} (\vartheta_{t},\X_{s^{-}}, \X_{s})}{ k_{ \gamma}( \vartheta,\X_{s^{-}}, \X_{s})}-1- \frac{ h^T}{  \sqrt{t}} \frac{ \nabla k_\gamma ( \vartheta,\X_{s^-},\X_s)}{ k_\gamma ( \vartheta,\X_{s^-},\X_s)}\right)dN^{\gamma}_{s}\\
+ \frac{ 1}{ \sqrt{t}}\int_{0} ^{t} \int_{E} h^T \nabla k_\gamma( \vartheta,\X_{s},\y)\nu_{ \gamma}(\X_{s},d\y) \gamma( \vartheta, \X_{s})ds,
\end{multline*}
\begin{multline*}
A_{t}=\frac{ 1}{ 2t}\int_{0} ^{t} \int_{E}\1_{k_\gamma>0} \frac{\left(h^T\nabla k_\gamma\right)^2}{ k_\gamma} ( \vartheta,\X_{s},\y) \nu_{ \gamma}(\X_{s},d\y) \gamma( \vartheta, \X_{s})ds\\
 - \frac{ 1}{ 2t}\int_{0}^{t} \frac{\left(h^T\nabla k_\gamma\right)^2}{ k_\gamma^2} ( \vartheta,\X_{s^-},\X_s)dN^{\gamma}_{s},
\end{multline*}
\begin{multline*}
R_{t}=\int_{0}^{t} \left(  \log\left( \frac{ k_{ \gamma}(\vartheta_{t},\X_{s^{-}}, \X_{s})}{ k_{ \gamma}( \vartheta,\X_{s^{-}}, \X_{s})}\right)-\frac{ k_{ \gamma}(\vartheta_{t},\X_{s^{-}}, \X_{s})}{ k_{ \gamma}( \vartheta,\X_{s^{-}}, \X_{s})}+1\right. \\ 
\left.+ \frac 1 {2t} \frac{\left(h^T\nabla k_\gamma\right)^2}{ k_\gamma^2} ( \vartheta,\X_{s^-},\X_s) \right) dN^{\gamma}_{s}.
\end{multline*}
We show that each of these three terms tends to 0 in $\P_{\x_0}^\vartheta$-probability as in the proof of Lemma~\ref{lem_gamma}. 

\smallskip

\noindent\underline{First step:} We have
$$[S_t]=\int_{0}^{t} \left( \frac{ k_{ \gamma} (\vartheta_{t},\X_{s^{-}}, \X_{s})}{ k_{ \gamma}( \vartheta,\X_{s^{-}}, \X_{s})}-1- \frac{ h^T}{  \sqrt{t}} \frac{ \nabla k_\gamma ( \vartheta,\X_{s^-},\X_s)}{ k_\gamma ( \vartheta,\X_{s^-},\X_s)}\right)^2dN^{\gamma}_{s}
$$
so that by Lemma~\ref{martingale_double} in Appendix~\ref{appendix martingale},
\begin{multline*}\langle S\rangle _{t}=\int_{0}^{t}  \int_{E}\left( \frac{ k_{ \gamma}(\vartheta_{t},\X_{s},\y)}{ k_{ \gamma}(\vartheta,\X_{s},\y)}-1-\frac{ h^T}{ \sqrt{t}}\frac{ \nabla k_\gamma ( \vartheta,\X_{s},\y)}{ k_\gamma ( \vartheta,\X_{s},\y)}\right)^2 \\\times k_{ \gamma} (\vartheta,\X_{s},\y) \nu_{ \gamma}(\X_{s},d\y)\gamma( \vartheta,\X_{s}) ds.\end{multline*}
Using the same inequalities as in \eqref{dl gamma} we obtain that
\begin{align*}
\langle S\rangle _{t}\leq \frac{ \|h\| ^2}{ t}\int_{0}^{t}   f^{k_{ \gamma}}_{ \vartheta}\left(\X_{s}, \frac{ \|h\|}{ \sqrt{t}}\right) \gamma(\vartheta, \X_{s})ds 
\end{align*}
where $f^{k_{ \gamma}}_{ \vartheta}$ is given by \eqref{f_kgamma}. We can then deduce that $\E_{\x_0}^\vartheta( \langle S\rangle _{t})$ tends to 0 by use of Conditions~\ref{cond:Btilde},  \ref{hypH} and the dominated convergence theorem, as in the proof of Lemma~\ref{lem_gamma}.

\smallskip

\noindent\underline{Second step:} 
The proof is identical to the second step of the proof of 
Lemma~\ref{lem_gamma}. It consists in splitting the integrals in $A_t$ according to whether $\left |\frac{h^T\nabla k_\gamma}{ k_\gamma} ( \vartheta,\X_{s^-},\X_s)\right|>\epsilon\sqrt t$
 or not. The convergence of $A_t$ is then proved similarly, by use of  Conditions~\ref{cond:Atilde}, \ref{hypH} and the dominated convergence theorem.

\smallskip

\noindent\underline{Third step:} The proof follows exactly the same lines as the proof of the third step 
of 
Lemma~\ref{lem_gamma} where $\gamma(.,.)$ is replaced by $k_{ \gamma}(.,.,.)$ and Condition~\ref{cond:Atilde} is used instead of Condition~\ref{cond:A}.

\subsection{Proof of Lemma~\ref{lem_L}}

Remember that $\X_t=\{(\Z^{(i)}_{t-T_i},\m^{(i)})\}$ for all $t\in[T_i,T_{i+1}[$, 
and that $L(\vartheta,\X_{[T_{i},T_{i+1}[})$ is a shortcut for $L(\Z^{(i)}_{[0,T_{i+1}-T_i]},\m^{(i)})$ given by \eqref{Lmove}. Accordingly, we have
\begin{multline*}\log L(\vartheta,\X_{[T_{i},T_{i+1}[})=\int_{T_i}^{T_{i+1}} \bar b_n(\vartheta,\X_s)^T a_{n}^{-1}( \X_s)d  \Z_{s-T_i}^{(i)}\\  
-\frac{1}{2}\int_{T_i}^{T_{i+1}}\bar b_n( \vartheta, \X_s)^T a_{n}^{-1}( \X_s) \bar b_n(\vartheta, \X_s)d s,
\end{multline*}
where $\bar b_n(\vartheta,\X_s)$ stands for $\bar b_n(\vartheta, \Z_{s-T_i}^{(i)},\m^{(i)})$, with $n=n(\X_{s})$,  and similarly for $a_{n}^{-1}( \X_s)$. 
Then we have
\begin{align*}
\log& \frac{L(\vartheta_t,\X_{[T_{i},T_{i+1}[})}{L(\vartheta,\X_{[T_{i},T_{i+1}[})} \\
=& \int_{T_i}^{T_{i+1}}( \bar b_n( \vartheta_t,\X_s)- \bar b_n(\vartheta,\X_s))^T a_{n}^{-1}( \X_s) (d  \Z_{s-T_i}^{(i)}  - \bar b_n( \vartheta, \X_s)ds) \\
&- \frac{1}{2}\int_{T_i}^{T_{i+1}}  ( \bar b_n( \vartheta_t,\X_s)- \bar b_n(\vartheta,\X_s))^T a_{n}^{-1}( \X_s)  ( \bar b_n( \vartheta_t,\X_s)- \bar b_n(\vartheta,\X_s) ds.
\end{align*}
Consequently, we may rewrite \eqref{lambdaGamma0} as
\begin{multline}\label{lambdaGamma}
\Lambda_{ L}^{ \vartheta_{t}/ \vartheta}(t)= \sum_{i\geq0}\int_{0}^t \1_{ [T_{i}; T_{i+1}[}(s)( \bar b_n( \vartheta_t,\X_s)- \bar b_n(\vartheta,\X_s))^T a_{n}^{-1}( \X_s)  dV^{(i)}_{s-T_i}\\
- \frac{1}{2}\int_{0}^{t}  ( \bar b_n( \vartheta_t,\X_s)- \bar b_n(\vartheta,\X_s))^T a_{n}^{-1}( \X_s)  ( \bar b_n( \vartheta_t,\X_s)- \bar b_n(\vartheta,\X_s)) ds.
\end{multline}
where $V_s^{(i)}=0$ if $s<0$ and 
\begin{align} \label{marting:Mk}
V^{(i)}_{s}= \Z^{(i)}_{s}-\Z^{(i)}_{0} -\displaystyle\int_{0}^{s} \bar b_n( \vartheta, \Z^{(i)}_u)du, \qquad s\geq 0.
\end{align}

We know  that $\Z^{(i)}$  is the solution of $M^{|n}(\z^{(i)},\m^{(i)})$ given by \eqref{EDSglob}, where $(\z^{(i)},\m^{(i)})$ are such that $\X_{T_i}=\{(\z^{(i)},\m^{(i)})\}$ and $n=n(\X_{T_i})$.
This implies that  $(V_s^{(i)})$ is a martingale with respect to $\mathcal F^{\Z^{(i)}}_t$, the natural filtration associated to $\Z^{(i)}$. Note that $(V_s^{(i)})$ is independent from $T_i$ and that $\mathcal F_t$ is generated by $\left\{T_{N_t},\X_{T_{N_t}}, \Z^{(N_t)}_{[0,t-T_{N_s}]}, (T_i,\X_{T_i},\Z^{(i)}_{[0,T_{i+1}-T_i]})_{i=1,\dots,N_t-1}\right\}$. This implies that for $T_i\leq t <T_{i+1}$, if $s\in [T_i,t]$,
$$\E(V^{(i)}_{t-T_i} |\mathcal F_s)=\E(V^{(i)}_{t-T_i} | T_i, \mathcal F^{\Z^{(i)}}_{s-T_i})=V^{(i)}_{s-T_i},$$
while if $s<T_i$, $\E(V^{(i)}_{t-T_i} |\mathcal F_s)=0$. 

Using these properties, we can easily verified that $M_{L}^{ \vartheta}$ given by \eqref{M_L}  is a $\mathcal F_{t}$-martingale. Similarly, since  $\langle V^{(i)}\rangle_s = \int_{0}^s a_n( \Z^{(i)}_u)du$, we have 
\begin{align*}
\langle M_{ L}^{ \vartheta}\rangle(t) &=\frac{1}{ t} \sum_{i\geq 0}\int_{0}^{t}\1_{ [T_{i}; T_{i+1}[}(s)  (\nabla \bar b_n( \vartheta,\X_{s}))^T a_{n}^{-1}( \X_{s})  \\
&\hspace{3cm}\times a_{n}( \X_{s}) a_{n}^{-1}( \X_{s})  ( \X_{s})\nabla \bar b_n( \vartheta,\X_{s}) ds\\
&=\frac{1}{ t}\displaystyle\int_{0}^{t}  (\nabla \bar b_n( \vartheta,\X_{s}))^T a_{n}^{-1}( \X_{s})\nabla \bar b_n( \vartheta,\X_{s}) ds.
\end{align*}

From this representation, \eqref{M_L} and \eqref{lambdaGamma}, we have
\begin{align*}
Rem^{L}_{ \vartheta,h}(t) = \Lambda_{L }^{ \vartheta_{t}/ \vartheta}(t) - h^T M_{ L}^{ \vartheta}(t) + \frac{ 1}{ 2}h^T\langle M_{ L}^{ \vartheta} \rangle(t) h
=S_t -\frac 12 R_t,\end{align*}
where
\begin{multline*}
S_t= \sum_{i\geq0}\int_{0}^{t} \1_{ [T_{i}; T_{i+1}[}(s)\bigg( \bar b_n( \vartheta_t,\X_s)- \bar b_n(\vartheta,\X_s)\\- \nabla \bar b_n( \vartheta,\X_{s})\frac{h}{\sqrt t}\bigg)^Ta_{n}^{-1}( \X_{s})dV^{(i)}_{s-T_i},
\end{multline*}
and
\begin{multline*}
R_t=\int_{0}^{t} \bigg( ( \bar b_n( \vartheta_t,\X_s)- \bar b_n(\vartheta,\X_s))^T a_{n}^{-1}( \X_s)  ( \bar b_n( \vartheta_t,\X_s)- \bar b_n(\vartheta,\X_s))  \\
 -\frac{1}{ t} h^T (\nabla \bar b_n( \vartheta,\X_{s}))^T a_{n}^{-1}( \X_{s})\nabla \bar b_n( \vartheta,\X_{s}) h\bigg) ds.
\end{multline*}

Let us prove that both $S_t$ and $R_t$ tend to 0 in $\P_{\x_0}^\vartheta$-probability. For $S_t$, using the Markov inequality and the fact that $\E_{\x_0}^\vartheta((S_t)^2)=\E_{\x_0}^\vartheta(\langle S_t\rangle)$, this boils down to  proving that $\E_{\x_0}^\vartheta(\langle S_t\rangle)$ tends to 0. 
We have, since $a_{n}= \bar\sigma_{n} \bar\sigma_{n}^T$, 
\begin{align*}
\langle S_t\rangle&= \int_{0}^{t}   \left( \bar b_n( \vartheta_t,\X_s)- \bar b_n(\vartheta,\X_s)- \nabla \bar b_n( \vartheta,\X_{s})\frac{h}{\sqrt t}\right)^Ta_{n}^{-1}( \X_{s})\\
&\hspace{3cm} \times\left( \bar b_n( \vartheta_t,\X_s)- \bar b_n(\vartheta,\X_s)- \nabla \bar b_n( \vartheta,\X_{s})\frac{h}{\sqrt t}\right) ds\\
&= \int_{0}^{t} \left\| \bar\sigma_n(\X_s)^{-1} \left( \bar b_n( \vartheta_t,\X_s)- \bar b_n(\vartheta,\X_s)- \nabla \bar b_n( \vartheta,\X_{s})\frac{h}{\sqrt t}\right) \right\|^2 ds.
\end{align*}
Using the same argument as for \eqref{dl gamma}, we get
 \begin{align*}
\langle S_t\rangle&\leq \frac{ \|h\| ^2}{ t} \int_{0}^{t} \left\| \bar\sigma_n(\X_s)^{-1}  \right\|^2 f^{\text{move}}_{ \vartheta}\left(\X_{s},\frac{ \|h\|}{ \sqrt{t}}\right) ds.
\end{align*}
where $f^{\text{move}}_{ \vartheta}$ is given by \eqref{f_L}. We then deduce that $\E_{\x_0}^\vartheta( \langle S\rangle _{t})$ tends to 0 thanks to  the dominated convergence theorem, using  Conditions~\ref{cond:Bmove} and \ref{hypH} as in the first step of the proof of Lemma~\ref{lem_gamma}. 

Concerning $R_t$, note that 
\begin{multline*}
R_t=\int_{0}^{t} \left\| \bar\sigma_n(\X_s)^{-1} \left( \bar b_n( \vartheta_t,\X_s)- \bar b_n(\vartheta,\X_s)\right)\right\|^2 \\ 
- \frac 1 t \left\| \bar\sigma_n(\X_s)^{-1} \nabla \bar b_n( \vartheta,\X_{s})h \right\|^2 ds,
\end{multline*}
that reads $\int(\|f\|^2-\|g\|^2)ds$ with obvious notation, so that using the same argument as in the  third step of the proof of Lemma~\ref{lem_gamma} (term $R_2$), it suffices to show that 
$\int \|f-g\|^2 ds =o(1)$ and $\int \|g\|^2ds=O(1)$ in $\P_{\x_0}^\vartheta$-probability. But  the former is exactly $\langle S_t\rangle$ studied above, that has been proven to converge to 0 in probability. And the expectation of the latter is bounded by use of Condition~\ref{cond:Amove}, proving that it is $O(1)$ in $\P_{\x_0}^\vartheta$-probability.  This concludes the proof.

\section{Ancillary results}
\subsection{Conditions for geometric ergodicity}\label{appendix:ergodicity}
The following proposition provides some conditions  ensuring the hypothesis~\ref{hypH},  namely non-explosion and geometric ergodicity of the BDMM process. It is based on the study conducted in \cite{emilien} for birth-death-move processes (without mutations). The arguments in presence of mutations are similar and sketched below. The conditions include the technical assumption that the process is Feller. This property is discussed in \cite{emilien} and proved for several examples of transition kernels that include all examples of this paper, provided the underlying spatial state $\Lambda$ is compact. The other conditions deal with the intensity functions of the jumps. First, the total intensity function $\alpha$ is assumed to be bounded, to avoid explosion of the process. Second, the death intensity function $\delta$  must compensate in a proper way the birth intensity function $\beta$. Notably, the ergodic properties do not depend on the inter-jump move process, neither on the mutation dynamics. 
We introduce the following notation:
\begin{equation} \label{defi_bn_dn}
    \beta_n = \underset{\x \in E_n}{\sup} \beta(\x), \quad \delta_n = \underset{\x \in E_n}{\inf} \delta(\x)\quad \textnormal{and}\quad \alpha_n = \beta_n + \delta_n.
\end{equation}
In fact the conditions of geometric ergodicity below are exactly the same as the conditions of ergodicity for a simple birth-death process $(\eta_t)_{t\geq 0}$ on $\N$ with birth rate $(\beta_n)$ and death rate $(\delta_n)$, as established in 
\cite{karlin1957}. This is because to address the ergodic properties of a BDMM $(\X_t)_{t\geq 0}$, we construct  a coupling between $(\X_t)_{t\geq 0}$ and $(\eta_t)_{t\geq 0}$, in such a way that $\eta_t=0$ implies $\X_t=\text{\O}$. Consequently $\text{\O}$ becomes a positive recurrent state for $(\X_t)_{t\geq 0}$ whenever 0 is a positive recurrent state for $(\eta_t)_{t\geq 0}$, which implies the geometric ergodicity of $(\X_t)_{t\geq 0}$.

\begin{prop} \label{prop:ergodicity}
	Let $(\X_t)_{t \geq 0}$ be a Feller BDMM process with a bounded intensity function $\alpha$.  Suppose that $\delta_n>0$ for all $n \geq 1 $ and one of the following condition holds:
	\begin{align} 
(i)&  \text{ there exists }  n_0 \geq 1 \text{ such that } \beta_{n}=0 \text{ for any } n \geq n_0,  \label{eq30}\\
(ii)& \ \beta_n>0 \text{ for all } n \geq 1, \ \sum_{n=2}^{\infty} \dfrac{\beta_1 \dots \beta_{n-1}}{\delta_1 \dots \delta_n} < \infty \text{ and }  \displaystyle \sum_{n=1}^{\infty} \dfrac{\delta_1  \dots \delta_n}{\beta_1  \dots \beta_n} = \infty,  \label{eq31}
\end{align}
where $\beta_n$ and $\delta_n$ are defined by \eqref{defi_bn_dn}. Assume  moreover that 
 \begin{align*}
 \sum_{n=2}^{\infty} \sqrt{\frac{\beta_1 \dots \beta_{n-1}}{\delta_1 \dots \delta_n}} < \infty \quad\text{ and }  \quad \exists N \geq 0, \; s.t. \; \forall \, n \geq N, \; \beta_n \leq \delta_{n+1}. 
 \end{align*}
	Then \ref{hypH} is satisfied.
\end{prop}

As mentioned above, following \cite{preston} and \cite{emilien}, the proof of this proposition is based on a coupling between  $(\X_t)_{t\geq 0}$ and $(\eta_t)_{t\geq 0}$, where the latter is a simple birth-death process with birth and death rates given by \eqref{defi_bn_dn}. We detail below how we construct the coupled process $\check C_t=(\X_t,\eta_t)$. This is a straightforward generalisation of the construction in \cite{emilien}, where we account for the presence of mutations. Specifically, $\check C$  is a jump move process on $E\times\N$ with intensity function $\check \alpha$ and transition kernel $\check K$ defined as follows. 
The intensity function $\ck \alpha : E \times \N \rightarrow [0,\infty)$ is given by 
\begin{equation*}
    \ck \alpha (\x,n)= \left \{    \begin{array}{lcl}
    \alpha(\x)+\alpha_n & & \textnormal{if } \x \in E_m, \, m \neq n  , \\
    \beta_n + \delta(\x)+\tau(\x) & & \textnormal{if } \x \in E_n.  
\end{array} \right.
\end{equation*}
Letting $K$ be the kernel given by \eqref{defK},  the transition kernel $\ck K$ takes the form, for any $A\subset E$:
\begin{enumerate}
    \item If $\x \in E_m, \, m \neq n$ : 
    \begin{align*}
        &\ck K((\x,n);A \times\{n\})  =  \dfrac{\alpha (\x)}{ \ck \alpha (\x,n)} K(\x,A);\\
        &\ck K((\x,n);\{\x\}\times\{n+1\})  =  \dfrac{\beta_n}{ \ck \alpha (\x,n)} ;\\
        &\ck K((\x,n);\{\x\}\times\{n-1\})  =  \dfrac{\delta_n}{ \ck \alpha (\x,n)}.\\
        \end{align*}
    \item If $\x \in E_n$ : 
   \begin{align*}
       &\ck K((\x,n);A\times\{n+1\})  =  \dfrac{\beta (\x)}{\ck \alpha (\x,n)} K_\beta(\x,A) ;\\
      &\ck K((\x,n);\{\x\}\times\{n+1\})  =  \dfrac{\beta_n-\beta(\x)}{\ck \alpha (\x,n)}; \\
       &\ck K((\x,n);A\times\{n-1\})  =  \dfrac{\delta_n}{\ck \alpha (\x,n)}K_\delta(\x,A); \\
      &\ck K((\x,n);A\times\{n\})  =  \dfrac{\delta(\x)-\delta_n}{\ck \alpha (\x,n)}K_\delta(\x,A)+ \dfrac{\tau(\x)}{\ck \alpha (\x,n)}K_\tau(\x,A) .
     \end{align*}
\end{enumerate}
The inter-jump move process of $\ck C$ is in turn a simple independent coupling between the move $\Y$ of $\X$ and a constant move on $\N$ (i.e. $y_t=y_0$ for all $t\geq 0$). 

The fact that the above construction is a proper coupling, in the sense that $\X$ and $\eta$ do constitute the marginal distributions of $\ck C$, can be proven exactly as in \cite{emilien} for the case without mutations. On the other hand, the key point is that if at some point, $C_t=(\X_t,\eta_t)$ is such that $n(\X_t)\leq \eta_t$, then by the above construction, $n(\X_s)\leq \eta_s$ for all $s\geq t$, almost surely. By this property, $\eta_t=0$ implies $\X_t=\text{\O}$ and the geometric ergodic conditions for $\eta$ stated in the proposition, coming from \cite{karlin1957},  are sufficient for the geometric ergodicity of $\X$. The rigorous proof of these claims can be found in \cite{emilien}.

We conclude this section by stating an integrated version of the geometric ergodic condition in  \ref{hypH}, which turns out to be useful in the proof of Theorem~\ref{LAN}. 
\begin{cor}\label{corH}
 Assume  \ref{hypH}. Then there exist a measure $\mu_{ \infty}^\vartheta$ on $E$,  and $c : E \to (0,\infty)$, possibly depending on $\vartheta$,  such that for any $t>0$, any $\x_0 \in E$ and any measurable and bounded function $g $ on $E$,
$$ \displaystyle\left| \frac{ 1}{ t} \E_{\x_0}^\vartheta \left(\int_{0}^tg(\X_{s})ds\right)-\int_{E}g(\y)\mu_{\infty}^\vartheta(d\y)\right|\leq \frac{ c(\x_0) \|g\|_{\infty}}{ t}.$$ 
\end{cor}

\begin{proof}
For any measurable and bounded function $g $ on $E$, using \ref{hypH},
\begin{align*}
\bigg| \frac{ 1}{ t} \E_{\x_0}^\vartheta \left(\int_{0}^tg(\X_{s})ds\right)-&\int_{E}g(\y)\mu_{\infty}^\vartheta(d\y)\bigg|
\\&\leq \frac{ 1}{ t} \int_0^t \left| \E_{\x_0}^\vartheta (g(\X_{s}) - \int_{E}g(\y)\mu_{\infty}^\vartheta(d\y)\right| ds\\ 
&\leq \frac{ c(\x_0)  \|g\|_{\infty}}{ t} \int_0^t r^s ds \leq  \frac{ c(\x_0)  \|g\|_{\infty}}{ t} \int_0^\infty r^s ds.
\end{align*}
Denoting (abusively) by  $c(\x_0)$ the constant $c(\x_0) \int_0^\infty r^s ds$, we get the result.
\end{proof}

\subsection{Some complements on martingales}\label{appendix martingale}

This section recalls some basic facts about martingales and their brackets. For general background and further details, we refer the reader to \cite{klebaner}.

The process $M_t$ is a martingale with respect to its history $\mathcal F_t$ if for any $s<t$, $\E(M_t|\mathcal F_s)=M_s$ almost surely. Its quadratic variations, also called square brackets, are defined by the limit in probability:
$$[M]_t = \lim \sum_{i=1}^n (M_{t_i^n}-M_{t_{i-1}^{n}})^2$$
where the limit is taken over partitions $0=t_0^n<t_1^n<\dots<t^n_n=t$ with $\max (t_i^n -  t_{i-1}^n)\to 0$, see \cite[Chapter 7.6]{klebaner}.
 For example, for a counting process $N_t$ with $\mathcal F_t$-intensity $\gamma$, the process $M_t=N_t- \int_0^t \gamma(s) ds$ is a $\mathcal F_t$-martingale with quadratic variations $[M]_t=N_t$. 
As another important example,  consider the stochastic integral $\int_0^t H(s) dM_s$, where $M_t$ is an  $\mathcal F_t$-martingale and the process $H_t$ is $\mathcal F_t$-predictable, see \cite[Chapter 8.4]{klebaner} for a definition. Then its quadratic variations are $$\left[\int_0^{\boldsymbol{\cdot}} H_s dM_s\right]_t=\int_0^t H_s^2 d[M]_s.$$

The square brackets $[M]_t$ are not necessarily predictable, as illustrated  by the above example involving the counting process $N_t$. By contrast, the angle (or sharp) brackets $\langle M\rangle_t$ of a martingale can be viewed as a predictable version of $[M]_t$. More precisely, the angle brackets are defined as the unique $\mathcal F_t$-predictable process $\langle M\rangle_t$ such that $[M]_t-\langle M\rangle_t$ is a $\mathcal F_t$-martingale; see \cite[Chapter 8.9]{klebaner}. If $M_t$ is continuous, $[M]_t=\langle M\rangle_t$, but the two notions differ in general. 
For instance, for the martingale $M_t=N_t- \int_0^t \gamma(s) ds$, where  $N_t$ is a counting process with intensity $\gamma$, we have $\langle M\rangle_t=\int_0^t \gamma(s) ds$. 

An important property, sometimes used as a definition, is the following: for any square integrable martingale, $\langle M\rangle_t$ is the unique predictable increasing process such that $M^2_t - \langle M\rangle_t$ is a martingale; see \cite[Theorem~8.24]{klebaner}.

We conclude this section by the following result, used in the proof of Theorem~\ref{LAN}. It is proven in Lemma~54 of \cite{these_Emilien}, the BDMM process being a particular case of a jump-move process studied in this work. 
 It is also stated under a slightly different setting in Proposition~3.3~(b) of \cite{eva2002}. 
\begin{lem}\label{martingale_double}
Let $(\X_t)_{t\geq 0}$ be a BDMM process on $E$ and denote its natural filtration by $(\mathcal F_t)_{t\geq 0}$. Let $g$ be a measurable bounded function defined on $E\times E$ and introduce for any $t>0$
$$M^*_t=\int_0^t g(\X_{s^-},\X_s)dN^\gamma_s - \int_0^t \gamma(\X_s) \int_E g(\X_s,\y)k_\gamma(\X_s,\y)\nu_\gamma(\X_s,d\y) ds.$$
If $N_t^\gamma<\infty$ for any $t\geq 0$, then $(M^*_t)_{t\geq0}$ is a $\mathcal F_t$-martingale. 
\end{lem}

\bibliographystyle{acm}
\bibliography{biblio}
\end{document}